\documentclass[11pt, letterpaper]{article}
\usepackage{times}
\usepackage[T1]{fontenc}
\usepackage{mathrsfs}
\usepackage{latexsym}
\usepackage[dvips]{graphics}
\usepackage{epsfig}
\usepackage{amsmath,amsfonts,amsthm,amssymb,amscd}
\input amssym.def
\input amssym.tex
\usepackage{color}
\usepackage[all]{xypic}

\newcommand\be{\begin{equation}}
\newcommand\ee{\end{equation}}
\newcommand\bea{\begin{eqnarray}}
\newcommand\eea{\end{eqnarray}}
\newcommand\bi{\begin{itemize}}
\newcommand\ei{\end{itemize}}
\newcommand\ben{\begin{enumerate}}
\newcommand\een{\end{enumerate}}
\newcommand\bc{\begin{center}}
\newcommand\ec{\end{center}}
\newcommand\ba{\begin{array}}
\newcommand\ea{\end{array}}

\newtheorem{thm}{Theorem}[section]

\newtheorem{cor}[thm]{Corollary}
\newtheorem{corollary}[thm]{Corollary}
\newtheorem{lem}[thm]{Lemma}
\newtheorem{lemma}[thm]{Lemma}
\newtheorem{prop}[thm]{Proposition}

\newtheorem{rek}[thm]{Remark}

\newtheorem{cla}[thm]{Claim}

\newcommand{\R}{\ensuremath{\mathbb{R}}}
\newcommand{\C}{\ensuremath{\mathbb{C}}}

\DeclareMathOperator{\End}{End}
\DeclareMathOperator{\Hom}{Hom}
\DeclareMathOperator{\ad}{ad}
\DeclareMathOperator{\YMH}{YMH}

\DeclareMathOperator{\QH}{QH}
\DeclareMathOperator{\re}{Re}
\DeclareMathOperator{\tr}{tr}
\DeclareMathOperator{\rank}{rank}
\DeclareMathOperator{\grad}{grad}
\DeclareMathOperator{\Lie}{Lie}
\DeclareMathOperator{\UT}{UT}
\DeclareMathOperator{\im}{im}
\DeclareMathOperator{\id}{id}
\DeclareMathOperator{\Gr}{Gr}
\DeclareMathOperator{\HNS}{HNS}

\DeclareMathOperator{\SU}{SU}
\DeclareMathOperator{\SL}{SL}
\DeclareMathOperator{\GL}{GL}
\DeclareMathOperator{\hyperquotient}{/ \negthinspace / \negthinspace /}

\newcommand{\muone}{F_A + [\phi, \phi^*]}

\newcommand{\starbar}{\mathop{\bar{*}}}
\newcommand{\mhiggs}{\mathcal{M}^{Higgs}}
\newcommand{\totalspace}{\mathcal{A} \times \Omega^{1,0}(\End(E))}
\newcommand{\holspace}{\mathcal{A}^{0,1} \times \Omega^{1,0}(\End(E))}
\newcommand{\linearisation}{\left. \frac{d}{dt} \right|_{t=0} }
\newcommand{\hone}[1]{\left\| #1 \right\|_{H^1}}
\newcommand{\hknorm}[1]{\left\| #1 \right\|_{H^k}}
\newcommand{\ltwo}[1]{\left\| #1 \right\|_{L^2}}

\begin{document}

\title{Morse Theory for the Space of Higgs Bundles}
\author{Graeme Wilkin}
\date{}



\maketitle

\begin{abstract}
The purpose of this paper is to prove the necessary analytic results to construct a Morse theory for the Yang-Mills-Higgs functional on the space of Higgs bundles over a compact Riemann surface. The main result is that the gradient flow with initial conditions $(A'', \phi)$ converges to a critical point of this functional, the isomorphism class of which is given by the graded object associated to the Harder-Narasimhan-Seshadri filtration of $(A'', \phi)$. In particular, the results of this paper show that the failure of hyperk\"ahler Kirwan surjectivity for rank $2$ fixed determinant Higgs bundles does not occur because of a failure of the existence of a Morse theory.
\end{abstract}

\section{Introduction}

This paper studies the convergence properties of the gradient flow of the Yang-Mills-Higgs functional on the space of Higgs bundles over a compact Riemann surface, as introduced by Hitchin in \cite{Hitchin87}. Higgs bundles that minimise this functional correspond to solutions of Hitchin's self-duality equations, which (modulo gauge transformations) correspond to points of the $\SL(n, \C)$ or $\GL(n, \C)$ character variety of the surface. The results of this paper provide the analytic background for the use of Morse theory in the spirit of Atiyah and Bott's approach for holomorphic bundles in \cite{AtiyahBott83} to compute topological invariants of these character varieties, a program that has been carried out for the case $n=2$ by the author, Georgios Daskalopoulos and Jonathan Weitsman in the paper \cite{dww}.

To precisely define the spaces and functions under consideration we use notation as follows. Let $X$ be a compact Riemann surface of genus $g$, and fix a $C^\infty$ complex vector bundle $E$ of rank $r$ and degree $d$ over $X$ with a Hermitian metric on the fibres. Let $\mathcal{A}$ denote the space of connections on $E$ compatible with the metric, and note that $\mathcal{A}$ is isomorphic to the space $\mathcal{A}^{0,1}$, the space of holomorphic structures on $E$. A pair $(A'', \phi) \in \mathcal{A}^{0,1} \times \Omega^{1,0}(\End(E)) \cong T^* \mathcal{A}$ is called a \emph{Higgs pair} if the relation $d_A'' \phi = 0$ is satisfied. Let $\mathcal{B}(r,d)$ denote the space of all Higgs pairs on $E$, this space can be visualised as follows. There is a projection map $p : \mathcal{B}(r,d) \rightarrow \mathcal{A}^{0,1}$ given by $p(A'', \phi) = A''$, the fibres of $p$ are vector spaces $\{ \phi \, | \, d_A'' \phi = 0 \}$, which change in dimension as the holomorphic structure changes. In this way it is easy to see that the space $\mathcal{B}(r,d)$ is singular. If the determinant of $E$ is held fixed throughout this process then the gauge group $\mathcal{G}$ has an $SU(r)$ structure, the space $\mathcal{A}$ consists of holomorphic structures with fixed determinant, and the Higgs field $\phi$ is also trace-free. This is known as the \emph{fixed determinant case}. If the determinant of $E$ is unrestricted then the gauge group $\mathcal{G}$ has a $U(r)$ structure and this is known as the \emph{non-fixed determinant case}.

In the following, $\mathcal{B}$ will be used to denote the space of Higgs bundles and the extra notation for the rank and degree of $E$ will be omitted if the meaning is clear from the context. $\mathcal{B}^{st}$ (resp. $\mathcal{B}^{ss}$) denotes the space of \emph{stable} (resp. \emph{semistable}) Higgs bundles, those for which every $\phi$-invariant holomorphic sub-bundle $F \subset E$ satisfies $\frac{\deg(F)}{\rank(F)} < \frac{\deg(E)}{\rank(E)}$ (resp. $\frac{\deg(F)}{\rank(F)} \leq \frac{\deg(E)}{\rank(E)}$). The \emph{moduli space of semistable Higgs bundles} is the space $\mhiggs(r, d) = \mathcal{B}^{ss} / \negthinspace / \mathcal{G}$, where the GIT quotient $/ \negthinspace /$ identifies the orbits whose closures intersect. In the fixed determinant case the moduli space is denoted $\mhiggs_0(r,d)$.

As noted in \cite{Hitchin87}, the space $T^* \mathcal{A}$ is an infinite-dimensional hyperk\"ahler manifold, and the action of the gauge group $\mathcal{G}$ induces three moment maps $\mu_1$, $\mu_2$ and $\mu_3$ taking values in $\Lie(\mathcal{G})^* \cong \Omega^2(\End(E))$ and given by
\begin{align*}
\mu_1(A, \phi) & = F_A + [\phi, \phi^*] \\
\mu_\C(A, \phi) & = \mu_2 + i \mu_3 = 2i d_A'' \phi
\end{align*}
A theorem of Hitchin in \cite{Hitchin87} and Simpson in \cite{Simpson88} identifies the moduli space of semistable Higgs bundles with the quotient $\left(\mu_1^{-1}(\alpha) \cap \mu_\C^{-1}(0)\right) / \mathcal{G}$, where $\alpha$ is a constant multiple of the identity that minimises $\| \mu_1 \|^2$, and which is determined by the degree of the bundle $E$. This is the \emph{hyperk\"ahler quotient} (as defined in \cite{HKLR87}) of $T^*\mathcal{A}$ by $\mathcal{G}$ at the point $(\alpha, 0, 0) \in \Lie(\mathcal{G})^* \otimes_\R \R^3$. 

The functional $\YMH(A, \phi) = \left\| \muone \right\|^2$ is defined on $\mathcal{B}$ using the $L^2$ inner product $\left< a, b \right> = \int_X \tr{a \starbar b}$. The purpose of this paper is to use the gradient flow of $\YMH$ to provide an analytic stratification of the space $\mathcal{B}$ for any rank and degree, and for both fixed and non-fixed determinant. The theorem of Hitchin and Simpson described above identifies the minimal stratum with the space of semistable Higgs bundles, the results here complete this picture by providing an algebraic description of the non-minimal strata for the flow in terms of the Harder-Narasimhan filtration.

\begin{thm}[Convergence of Gradient Flow]\label{thm:convergenceintro}
The gradient flow of $$\YMH(A, \phi) = \left\| F_A + [\phi, \phi^*] \right\|^2$$ converges in the $C^\infty$ topology to a critical point of $\YMH$. Moreover, let $r(A_0, \phi_0)$ be the map which take the initial conditions $(A_0, \phi_0)$ to their limit under the gradient flow equations. Then for each connected component $\eta$ of the set of critical points of $\YMH$, the map $r:  \{(A_0, \phi_0) \in \mathcal{B} : r(A_0, \phi_0) \in \eta  \} \rightarrow \eta$ is a $\mathcal{G}$-equivariant continuous map.
\end{thm}

This theorem is proved in Section \ref{sec:analysis}. On each non-minimal critical set, the critical point equations of $\YMH$ define a splitting of $E = F_1 \oplus \cdots \oplus F_n$ into $\phi$-invariant holomorphic sub-bundles. The degree of each component of the splitting is (up to re-ordering) well-defined on each connected component of the set of critical points, and each component can be classified by the Harder-Narasimhan type of the splitting into sub-bundles. This leads to the following stratification of the space $\mathcal{B}$.

\begin{cor}[Description of Analytic Stratification]\label{cor:stratificationintro}
The space $\mathcal{B}$ admits a stratification in the sense of \cite{AtiyahBott83} Proposition 1.19 $(1)$-$(4)$, which is indexed by the set of connected components of the critical points of the functional $\YMH$. 
\end{cor}
As described in \cite{HauselThaddeus04}, $\mathcal{B}$ can also be stratified algebraically by the $\phi$-invariant Harder-Narasimhan type of each Higgs bundle. The following theorem shows that this stratification is the same as that in Corollary \ref{cor:stratificationintro}.

\begin{thm}[Equivalence of Algebraic and Analytic Stratifications]\label{thm:equivalenceintro}
The algebraic stratification of $\mathcal{B}$ by Harder-Narasimhan type is equivalent to the analytic stratification of $\mathcal{B}$ by the gradient flow of the functional $\YMH$.
\end{thm}
This theorem is proved in Section \ref{sec:stratification}. Moreover, the following theorem (proved in Section \ref{sec:gradedobject}) provides an algebraic description of the limit of the gradient flow in terms of the Harder-Narasimhan-Seshadri filtration of the bundle.
\begin{thm}[Convergence to the graded object of the HNS filtration]\label{thm:gradedobjectintro}
The isomorphism class of the retraction $r: \mathcal{B} \rightarrow \mathcal{B}_{crit}$ onto the critical sets of $\YMH$ is given by
\begin{equation}
r(A'', \phi) \cong \Gr^{\HNS}(A'', \phi)
\end{equation}
where $\Gr^{\HNS}(A'', \phi)$ is defined in Section \ref{sec:gradedobject}.
\end{thm}

A long-standing question for finite-dimensional hyperk\"ahler quotients $M \hyperquotient G$ is the question of whether the hyperk\"ahler Kirwan map is surjective. In infinite dimensions this is not true, since a comparison of the Betti numbers from the computation of $P_t(\mhiggs_0(2,1))$ in \cite{Hitchin87}, together with the calculation of $P_t(B \mathcal{G}^{\SU(2)})$ from Theorem 2.15 of \cite{AtiyahBott83}, shows that the hyperk\"ahler Kirwan map $\kappa_{HK} : H_\mathcal{G}^*(T^* \mathcal{A}) \rightarrow H_\mathcal{G}^*(\mu_1^{-1}(\alpha) \cap \mu_\C^{-1}(0))$ cannot be surjective in the case of rank $2$ degree $1$ fixed determinant Higgs bundles. It would have been reasonable to conjecture that this failure of surjectivity occurs because of a failure of the Morse theory for this infinite-dimensional example, however the results of this paper show that the Morse theory actually does work, and the paper \cite{dww} explains the failure of hyperk\"ahler Kirwan surjectivity for this example in terms of the singularities in the space $\mathcal{B}$. 

The proof of Theorem \ref{thm:convergenceintro} is an extension of the approach of Rade in \cite{radethesis} and \cite{Rade92} where it was shown that the gradient flow of the Yang-Mills functional converges in the $H^1$ norm when $X$ is a $2$ or $3$ dimensional manifold, thus providing a purely analytic stratification of the space $\mathcal{A}$. Rade's proof was based on a technique of Simon in \cite{Simon83}, the key step being to show that a Lojasiewicz-type inequality holds in a neighbourhood of each critical point. Theorem \ref{thm:convergenceintro} extends this result to Higgs bundles and also improves on the convergence (showing $C^\infty$ convergence instead of $H^1$ convergence), by using a Moser iteration argument. 

This paper is organised as follows. Section \ref{sec:symplectic} sets the notation that is used in the rest of the paper. In Section \ref{sec:analysis} we prove the convergence result, Theorem \ref{thm:convergenceintro}. Section \ref{sec:stratification} contains the proof of the equivalence between the analytic stratification defined by the gradient flow of $\YMH$ and the algebraic stratification by Harder-Narasimhan type (Theorem \ref{thm:equivalenceintro}) and Section \ref{sec:gradedobject} shows that the gradient flow converges to the graded object of the Harder-Narasimhan-Seshadri double filtration (Theorem \ref{thm:gradedobjectintro}).

{\bf Acknowledgements:} I am indebted to my advisor Georgios Daskalopoulos for his advice and encouragement during the writing of this paper, and also to Jonathan Weitsman and Tom Goodwillie for many useful discussions. I would also like to thank the American Institute of Mathematics and the Banff International Research Station for their hospitality during the respective workshops "Moment Maps and Surjectivity in Various Geometries" in August 2004 and "Moment Maps in Various Geometries" in June 2005.

\section{Symplectic Preliminaries}\label{sec:symplectic}

In this section we derive the basic symplectic formulas that are used to set the notation and sign conventions for the rest of the paper. First identify 
\begin{equation*}
\totalspace \cong \holspace
\end{equation*}
where $\mathcal{A}^{0, 1}$ denotes the space of holomorphic structures on $E$ (as in \cite{AtiyahBott83} Section 5), and note that the tangent space is isomorphic to
\begin{equation}\label{eqn:holtangent}
T_{(A'', \phi)} \left( \holspace \right) \cong \Omega^{0, 1} (\End(E)) \times \Omega^{1, 0} (\End(E))
\end{equation}
The metric used here is given by
\begin{equation}\label{eqn:holmetric}
g\left( \left( \begin{matrix} a_1'' \\ \varphi_1 \end{matrix} \right), \left( \begin{matrix} a_2'' \\ \varphi_2 \end{matrix} \right) \right) = 2 \re \int_X \tr \{ a_1'' \starbar a_2'' \} + 2 \re \int_X \tr \{\varphi_1 \starbar \varphi_2 \}
\end{equation}
where $\starbar(\cdot) = *(\cdot)^*$, $*$ being the usual Hodge star operator and $(\cdot)^*$ the Hermitian adjoint with respect to the Hermitian metric on the fibres. Similarly, the inner product on $\Lie (\mathcal{G})$ is defined as follows
\begin{equation}\label{eqn:innerproduct}
\left< u, v \right> = \int_X \tr \{u \starbar v \} = -\int_X \tr \{u*v\}
\end{equation}
The dual pairing $\Lie(\mathcal{G})^* \times \Lie(\mathcal{G}) \rightarrow \R$ is given by
\begin{equation}\label{eqn:dualpairing}
\mu \cdot u = -\int_X \tr \{u \mu \}
\end{equation}
and noting that $\mu \cdot u = < u, *\mu >$ we see that the identification of $\Lie(\mathcal{G})^*$ with $\Lie(\mathcal{G})$ for this choice of inner product and dual pairing is the Hodge star operator $*:\Omega^2(\End(E)) \rightarrow \Omega^0(\End(E))$. The group action of $\mathcal{G}$ on $\holspace$ is given by
\begin{equation}\label{eqn:groupaction}
g \cdot \left( \begin{matrix} A'' \\ \phi \end{matrix} \right) = \left( \begin{matrix} g^{-1} A'' g + g^{-1} dg \\ g^{-1} \phi g \end{matrix} \right)
\end{equation}
Differentiating this gives us the infinitesimal action
\begin{equation}\label{eqn:infaction}
\rho_{(A'', \phi)}(u) = \left( \begin{matrix} d_A''u \\ [\phi, u] \end{matrix} \right) 
\end{equation}
The extra notation denoting the point $(A'', \phi)$ will be omitted if the meaning is clear from the context. If $\rho(u)=0$ then differentiating again gives us the infinitesimal action of $u$ on the tangent space $T_{(A'', \phi)} \left( T^* \mathcal{A}^{0,1} \right)$
\begin{equation}\label{eqn:inftangaction}
\delta \rho(u) \left( \begin{matrix} a'' \\ \varphi \end{matrix} \right) = \linearisation \left( \begin{matrix} d_{A+ta}'' u \\ [\phi + t \varphi, u] \end{matrix} \right) = \left( \begin{matrix} [a'', u] \\ [\varphi, u] \end{matrix} \right)
\end{equation}
For some calculations (such as those in Section \ref{sec:analysis}) it is more convenient to use the identification
\begin{align*}
T_{(A, \phi)} \left( \totalspace \right) & \cong \Omega^1(\End(E)) \\
\left( \begin{matrix} a \\ \varphi \end{matrix} \right) \mapsto & a + \varphi + \varphi^*
\end{align*}
where $a \in \Omega^1(\ad(E))$ and $\varphi \in \Omega^{1, 0}(\End(E))$. This allows us to consider a Higgs pair $(A, \phi)$ as a $GL(n, \C)$ connection on $E$, given by
\begin{align*}\label{eqn:glncconnection}
D_{(A, \phi)} : \Omega^0(\End(E)) & \rightarrow \Omega^1(\End(E)) \\
 u & \mapsto d_A u + [\phi + \phi^*, u]
\end{align*}
Note that if $u \in \Omega^1(\ad(E))$, then $D_{(A, \phi)} u - \left( D_{(A, \phi)} u \right)^* = 2d_A u$ and $D_{(A, \phi)} u + \left( D_{(A, \phi)} u \right)^* = 2[\phi + \phi^*, u]$, and therefore by splitting the tangent space into skew-adjoint and self-adjoint parts we can use this interpretation to give us the infinitesimal action of $\mathcal{G}$ on $\mathcal{A}^{GL(n, \C)}$, the space of $GL(n, \C)$ connections on $E$.
\begin{align*}
\rho: \Lie(\mathcal{G}) & \rightarrow \Omega^1(\End(E)) \\
\rho(u) & = D_{(A, \phi)} u
\end{align*}
In the case of a Higgs pair $(A, \phi)$ a simple computation shows that the curvature of $D_{(A, \phi)}$, denoted $F_{(A, \phi)}$, satisfies $F_{(A, \phi)} = F_A + [\phi, \phi^*]$. It is useful to note that $F_{(A, \phi)}^* = -F_{(A, \phi)}$. Now consider a general hyperk\"ahler manifold $M$ with the hyperhamiltonian action of a Lie group $G$. Let $\rho: \Lie(G) \rightarrow C^\infty(T M)$ be the infinitesimal action of $G$, and define $\rho_x^*$ to be the operator adjoint of $\rho_x$ at the point $x \in M$ with respect to the metric $g$ and the pairing $< \cdot, \cdot >$ on the space $\Lie(G)$
\begin{equation*}
g( \rho_x(u), X) = \left< u, \rho_x^*(X) \right>
\end{equation*}
The moment map condition $d{\mu_1}(X) \cdot u = \omega(\rho_x(u), X) = g(I\rho_x(u), X)$ shows that $d{\mu_1} \in \Lie(G)^*$ can be identified with $-\rho_x^*I \in \Lie(G)$. By differentiating the condition $\mu(g \cdot x) = g^{-1} \mu(x) g$ we obtain the following formula
\begin{equation}\label{eqn:infadjoint}
\rho_x^* I \rho_x(u) = -[* \mu_1(x), u]
\end{equation}
and similarly $\rho_x^* J \rho_x(u) = -[* \mu_2(x), u]$ and $\rho_x^* K \rho_x(u) = -[* \mu_3(x), u]$, where $*$ is used to denote the identification of $\Lie(G)$ with $\Lie(G)^*$. Differentiating again, we obtain the following product formulas for $\rho_x^*$ acting on $I \delta \rho_x(u)(X)$ and $\delta \rho_x(u)(X)$.
\begin{equation}\label{eqn:Iproductform}
\rho_x^* I \delta \rho_x(u)(X) = [\rho_x^*(IX), u] - (\delta \rho)_x^*\left(X, I \rho_x(u) \right)
\end{equation}
For the space of Higgs bundles with the action of $\mathcal{G}$ on the space $T^* \mathcal{A}$, a calculation shows that the complex structure $I$ commutes with the infinitesimal action on the tangent space in the following sense
\begin{equation}\label{eqn:IcommutesD}
I \delta \rho_x(u)(X) = \delta \rho_x(u)(IX)
\end{equation}
Therefore we can use \eqref{eqn:Iproductform} to derive the product formula
\begin{equation}\label{eqn:productform}
\rho_x^*  \delta \rho_x(u)(X) = [\rho_x^* X, u] + (\delta \rho)_x^*(X, \rho_x(u))
\end{equation}
Note that this formula is true for any K\"ahler manifold for which the commutativity relation \eqref{eqn:IcommutesD} holds.

\section{Convergence of the gradient flow}\label{sec:analysis}

Using the notation and formulae of the previous section, a calculation shows that for a K\"ahler manifold $M$ with moment map $\mu_1$ associated to a Hamiltonian $G$-action, the downwards gradient flow equations for the functional $\frac{1}{2} \left\| \mu_1(x) \right\|^2$ are given by $\frac{\partial x}{\partial t} = -I \rho_x(*\mu_1(x))$. More explicitly, for the functional $\YMH$ on the manifold $T^* \mathcal{A}^{0,1}$, the gradient flow equations are
\begin{align}\label{eqn:gradfloweqns}
\begin{split}
\frac{\partial A''}{\partial t} & = i d_A''*(F_A + [\phi, \phi^*]) \\ 
\frac{\partial \phi}{\partial t} & = i[\phi, *(F_A + [\phi, \phi^*])]
\end{split}
\end{align}
The purpose of this section is to prove the following theorem.
\begin{thm}[Convergence of Gradient Flow]\label{thm:convergence}
The gradient flow of $$\YMH(A, \phi) = \left\| F_A + [\phi, \phi^*] \right\|^2$$ with initial conditions in $\mathcal{B}$ converges in the $C^\infty$ topology to a critical point of $\YMH$. Moreover, let $r(A_0, \phi_0)$ be the map which takes the initial conditions $(A_0, \phi_0)$ to their limit under the gradient flow equations. Then for each connected component $\eta$ of the set of critical points of $\YMH$, the map $r:  \{(A_0, \phi_0) \in \mathcal{B} : r(A_0, \phi_0) \in \eta  \} \rightarrow \mathcal{B}$ is a $\mathcal{G}$-equivariant continuous map.
\end{thm}

In \cite{radethesis} and \cite{Rade92}, Rade proves convergence of the gradient flow of the Yang-Mills functional in the $H^1$ norm when the base manifold is $2$ or $3$ dimensional. Here we extend Rade's results to the case of Higgs bundles over a compact Riemann surface, and use a Moser iteration method to improve the regularity to smooth convergence. This relies on the following propositions.

\begin{prop}[Existence and Uniqueness]\label{lem:existence}
The gradient flow equations for the functional $\YMH$ have a unique solution which exists for all time. 
\end{prop}

\begin{prop}[Convergence modulo gauge transformations]\label{lem:convmodgauge}
For each $k>0$ there exist sequences $\{t_n\} \subseteq \R_+$ and $\{g_n\} \subseteq \mathcal{G}$ of Sobolev class $H^{k+2}$ such that $t_n \rightarrow \infty$ and $g_n \cdot \left(A(t_n), \phi(t_n) \right)$ converges strongly in the $H^k$ norm to a critical point $(A_\infty, \phi_\infty)$ of the functional $\YMH(A, \phi)$.
\end{prop}

\begin{prop}[Continuous dependence on initial conditions]\label{lem:continuousdependence}
For all $k \geq 1$ and $T > 0$, a solution to the gradient flow equations \eqref{eqn:gradfloweqns} at time $T$ depends continuously on the initial conditions in the topology induced by the $H^k$ norm.
\end{prop}

\begin{prop}[Lojasiewicz inequality]\label{lem:loja}
Given a critical point $(A_\infty, \phi_\infty)$ of the functional $\YMH$, there exists $\varepsilon_1 > 0$ such that the inequality
\begin{equation}\label{eqn:loja}
\| D_{(A, \phi)}^* F_{(A, \phi)} \|_{L^2} \geq c \left| \YMH(A, \phi) - \YMH(A_\infty, \phi_\infty) \right|^{1-\theta}
\end{equation}
holds for some $\theta \in (0, \frac{1}{2})$ whenever $\hone{(A, \phi) - (A_\infty, \phi_\infty)} < \varepsilon_1$.
\end{prop}
                                                                                                
\begin{prop}[Interior Estimate]\label{lem:energy}
Let $\varepsilon_1$ be as in Proposition \ref{lem:loja}, $k$ any positive integer and $S$ any real number greater than $1$. Given a critical point $(A_\infty, \phi_\infty)$ of the functional $\YMH$ and some $T$ such that $0 \leq T \leq S-1$, there exists a constant $c$ such that for any solution $(A(t), \phi(t))$ to the gradient flow of $\YMH(A, \phi)$ satisfying $\hknorm{(A(t), \phi(t))-(A_\infty,\phi_\infty)} < \varepsilon_1$ for all $t \in [T, S]$ then
\begin{equation}\label{eqn:energy}
\int_{T+1}^S \hknorm{\left( \frac{\partial A}{\partial t}, \frac{\partial \phi}{\partial t} \right)}\, dt \leq c\int_T^S \ltwo{\left( \frac{\partial A}{\partial t}, \frac{\partial \phi}{\partial t} \right)} \, dt
\end{equation}
\end{prop}

Assuming the results of these propositions, the proof of Theorem \ref{thm:convergence} proceeds as follows.

\begin{prop}\label{prop:convergence}
Let $(A_\infty, \phi_\infty)$ be a critical point of the functional $\YMH$ and let $k > 0$. Then there exists $\varepsilon > 0$ such that if $(A(t), \phi(t))$ is a solution to \eqref{eqn:gradfloweqns} and if for some $T \geq 0$
\begin{equation}\label{eqn:initallyclose}
\hknorm{(A(T), \phi(T)) - (A_\infty, \phi_\infty)} < \varepsilon
\end{equation}
then either $\YMH(A(t), \phi(t)) < \YMH(A_\infty, \phi_\infty)$ for some $t > T$, or $(A(t), \phi(t))$ converges in $H^k$ to a critical point $(A_\infty', \phi_\infty')$ as $t\rightarrow \infty$, where $\YMH(A_\infty', \phi_\infty') = \YMH(A_\infty, \phi_\infty)$. In the second case the following inequality holds
\begin{equation}\label{eqn:critpointsclose}
\hknorm{(A_\infty', \phi_\infty') - (A_\infty, \phi_\infty) } \leq c \hknorm{ (A(T), \phi(T)) - (A_\infty, \phi_\infty) }^{2\theta}
\end{equation}
with $\theta$ as in Proposition \ref{lem:loja} and where $c$ depends on the choice of critical point $(A_\infty, \phi_\infty)$.
\end{prop}

The method of proof of Proposition \ref{prop:convergence} is the same as the proof of Proposition 7.4 in \cite{Rade92}, and so it is omitted. Here we use Higgs bundles instead of connections, and also derive estimates in the $H^k$ norm using Proposition \ref{lem:energy}.

Using the above results we can now prove the main theorem of this section.

\begin{proof}[Proof of Theorem \ref{thm:convergence}]

Let $(A(t), \phi(t))$ be a solution to the gradient flow equations, and let $\mathcal{G}_{H^{k+2}}$ denote the completion of the group $\mathcal{G}$ in the $H^{k+2}$ norm. Proposition \ref{lem:convmodgauge} shows that there exists a sequence $\{t_n \}$ such that $t_n \rightarrow \infty$ and $\{g_n \} \subset \mathcal{G}_{H^{k+2}}$ such that 
\begin{equation}\label{eqn:convmodgauge}
g_n \cdot (A(t_n), \phi(t_n)) \rightarrow (A_\infty^k, \phi_\infty^k)
\end{equation}
strongly in $H^k$, where $(A_\infty^k, \phi_\infty^k)$ is a critical point of the functional $\YMH$. Since the functional $\YMH$ is invariant under the action of $\mathcal{G}$ and decreasing along the gradient flow then
\begin{equation}
\YMH(g_n \cdot (A(t_n), \phi(t_n))) \geq \cdots \geq \YMH(A_\infty^k, \phi_\infty^k)
\end{equation}
In particular $\YMH(A(t), \phi(t)) \geq \YMH(A_\infty^k, \phi_\infty^k)$ for all $t$. Equation \eqref{eqn:convmodgauge} implies that given any $\varepsilon$ there exists some $n$ such that 
$$\hknorm{g_n \cdot (A(t_n), \phi(t_n)) - (A_\infty^k, \phi_\infty^k)} < \varepsilon$$
The gradient flow equations are both unitary gauge-invariant and translation invariant with respect to $t$, and so $g_n \cdot (A(t_n + t), \phi(t_n + t))$ is also a solution. For notation let $(A'(t), \phi'(t)) = g_n \cdot (A(t_n + t), \phi(t_n + t))$. Then 
$$\hknorm{(A'(t), \phi'(t)) - (A_\infty^k, \phi_\infty^k)} < \varepsilon$$
for all $t \geq 0$, and $\YMH(A'(t), \phi'(t)) \geq \YMH(A_\infty^k, \phi_\infty^k)$. Therefore we are in the second case of Proposition \ref{prop:convergence}, and so $(A'(t), \phi'(t)) \rightarrow (A_\infty', \phi_\infty')$ strongly in $H^k$ for some critical point $(A_\infty', \phi_\infty')$. Therefore
\begin{eqnarray*}
g_n \cdot (A(t_n+t), \phi(t_n+t)) & \rightarrow & (A_\infty', \phi_\infty') \\
\Leftrightarrow \, \, (A(t), \phi(t)) & \rightarrow & (g_n)^{-1} \cdot (A_\infty', \phi_\infty')
\end{eqnarray*}
Since the critical point equations are $\mathcal{G}$ invariant, then $(g_n)^{-1} \cdot (A_\infty', \phi_\infty')$ is a critical point of the functional $\YMH$.

Therefore the gradient flow converges in $H^k$ to a critical point $(A_\infty^k, \phi_\infty^k)$ for all $k>0$. Since $\| \cdot \|_{H^k} \leq \| \cdot \|_{H^{k+1}}$ for all $k$ then $(A_\infty^k, \phi_\infty^k) = (A_\infty^{k+1}, \phi_\infty^{k+1}) = \cdots = (A_\infty, \phi_\infty)$ for all $k$. The Sobolev embedding theorem implies $C^{k-2} \subset H^k$ for all $k$, and so the gradient flow of $\YMH$ converges smoothly to $(A_\infty, \phi_\infty)$.
                                                                                                
To show that the limit depends continuously on the initial data, consider a solution $(A(t), \phi(t))$ to the gradient flow equations that converges in $H^k$ to a critical point $(A_\infty, \phi_\infty)$.   Since $(A(t), \phi(t))$ converges to $(A_\infty, \phi_\infty)$ then there exists $T$ such that $\hknorm{(A(T), \phi(T)) - (A_\infty, \phi_\infty)} < \frac{1}{2} \beta_2$.  Proposition \ref{lem:continuousdependence} states that finite time solutions to the gradient flow equations depend continuously on the initial conditions, therefore given $\beta_2$ and $T$ as above there exists $\beta_3>0$ such that if $\hknorm{(A'(0), \phi'(0)) - (A(0), \phi(0))}< \beta_3$ then 
\begin{equation*}
\hknorm{(A'(T), \phi'(T)) - (A(T), \phi(T))} < \frac{1}{2}\beta_2
\end{equation*}

It then follows from Proposition \ref{prop:convergence} that for any $\beta_1 > 0$ there exists $\beta_2 > 0$ such that if $(A'(t), \phi'(t))$ is another solution to the gradient flow equations which satisfies 
$$\hknorm{(A'(T), \phi'(T)) - (A_\infty, \phi_\infty)} < \beta_2$$
for some $T$, and which converges to $(A_\infty', \phi_\infty')$ in the same connected component of the set of critical points of $\YMH$ as $(A_\infty, \phi_\infty)$, then we have the estimate $\hknorm{(A_\infty', \phi_\infty') - (A_\infty, \phi_\infty)} < \beta_1$. Therefore, given any initial condition $(A(0), \phi(0))$, the above results show that for any $\beta_1>0$ there exists $\beta_3>0$ such that given another initial condition $(A'(0), \phi'(0))$ satisfying both 
$$\hknorm{(A'(0), \phi'(0)) - (A(0), \phi(0))}< \beta_3$$
and also that $r(A'(0), \phi'(0))$ and $r(A(0), \phi(0))$ are in the same connected component of the set of critical points of $\YMH$, then $(A'(t), \phi'(t))$ converges in $H^k$ to a critical point $(A_\infty', \phi_\infty')$ such that 
\begin{equation*}
\hknorm{(A_\infty', \phi_\infty') - (A_\infty, \phi_\infty)} < \beta_1 \qedhere
\end{equation*}

\end{proof}

\subsection{Existence and uniqueness of the gradient flow}\label{subsec:existence}

In this section we prove Proposition \ref{lem:existence}, which states existence and uniqueness for the gradient flow equations \eqref{eqn:gradfloweqns} with initial conditions $(A_0, \phi_0) \in \mathcal{B}$.

In \cite{Simpson88} the gradient flow equations \eqref{eqn:gradfloweqns} are studied as evolution equations on the space of Hermitian metrics on $E$. This equivalence is described as follows: fix a holomorphic structure on $E$ and a holomorphic section $\phi_0$ of $\Omega^{1,0}(\End(E))$. Now let $H$ be any Hermitian metric on $E$ and let $D_H = d'' + d_H' + \phi_0 + \phi_0^{*_H}$ be a $GL(r, \C)$ connection, where $d'' + d_H' = d_A$ denotes the metric connection on $E$ and $\phi_0^{*_H}$ is defined using Hermitian transpose with respect to the metric $H$. More explicitly, we can write
\begin{equation}\label{eqn:metricconnection}
d_A =  d'' + d' + H^{-1} d' H
\end{equation}
\begin{equation}\label{eqn:metricphi}
\phi_0 + \phi_0^{*_H} = \phi_0 + H^{-1} \bar{\phi_0}^T H
\end{equation}
Denote the curvature of $D_H$ by $F_H$ and let $\Lambda F_H^\perp = \Lambda F_H - \lambda \cdot \id$ where $\lambda =  \tr\{ F_H \}$ is a function $\lambda: X \rightarrow \C$, and $\Lambda : \Omega^k \rightarrow \Omega^{k-2}$ is defined in the standard way using the K\"ahler structure on $X$. For $X$ a compact Riemann surface, the following theorem is a special case of that given by Simpson in Section 6 of \cite{Simpson88}.

\begin{thm}[Simpson]\label{thm:simpsonexistence}
Solutions to the nonlinear heat equation
\begin{equation}\label{eqn:simpsonheat}
H^{-1} \frac{\partial H}{\partial t} = -2i \Lambda F_H^\perp = -2 i \left(\Lambda F_H - \lambda \cdot \id \right)
\end{equation}
exist for all time and depend continuously on the initial condition $H(0)$.
\end{thm}

The proof of Proposition \ref{lem:existence} relies on showing that equation \eqref{eqn:simpsonheat} is equivalent to the gradient flow of $\YMH$. As an intermediate step we use the following \emph{equivalent flow equations} for $(\tilde{A}(t), \tilde{\phi}(t))$
\begin{align}\label{eqn:equivflow}
\begin{split}
\frac{\partial \tilde{A''}}{\partial t} & = i\tilde{d_A}'' * (\tilde{F_A} + [\tilde{\phi}, \tilde{\phi}^*]) + \tilde{d_A}'' \alpha \\
\frac{\partial \tilde{\phi}}{\partial t} & = i[\tilde{\phi}, *(\tilde{F_A} + [\tilde{\phi}, \tilde{\phi}^*])] + [\tilde{\phi}, \alpha] 
\end{split}
\end{align}
for some one-parameter family $\alpha(t) \in \Omega^0(\ad(E))$. Note that the new terms in the equations correspond to the infinitesimal action of $\alpha$ at $(\tilde{A}'', \tilde{\phi})$. These equations are Higgs bundle versions of the equivalent flow equations used in \cite{Donaldson85} to prove existence for the Yang-Mills gradient flow equation, however here we also use the methods of \cite{Hong01} to show the relationship between the equivalent flow equations and the gradient flow equations. To achieve this let $H(t)= H_0 h(t)$, note that $h^{-1}\frac{\partial h}{\partial t} = H^{-1} \frac{\partial H}{\partial t}$ and consider the following equation for $h(t)$
\begin{equation}\label{eqn:gaugeflow}
\frac{\partial h}{\partial t} = -2i h *\left(F_{A_0} + d_{A_0}''(h^{-1}(d_{A_0}'h)) \\ + [\phi_0, h^{-1} \phi_0^* h] \right) + 2i\lambda h
\end{equation}
where $d_{A_0}$ is the metric connection for $H(0)$. The proof of Proposition \ref{lem:existence} requires the following lemmas, which together show that Theorem \ref{thm:simpsonexistence} implies existence for equation \eqref{eqn:equivflow}.
\begin{lem}\label{lem:gaugeflowtosimpson}
Existence for equation \eqref{eqn:simpsonheat} implies existence for equation \eqref{eqn:gaugeflow}.
\end{lem}

\begin{proof}[Proof of Lemma \ref{lem:gaugeflowtosimpson}]

By explicit computation using \eqref{eqn:metricconnection} and \eqref{eqn:metricphi} we also have
\begin{equation}\label{eqn:metriccurvature}
F_{H(t)} = F_{A_0} + d_{A_0}''(h^{-1}(d_{A_0}'h)) + [\phi_0, h^{-1} \phi_0^* h] \qedhere
\end{equation}
\end{proof}

Note that $h(0) = \id$ and that $h(t)$ is positive definite, therefore we can choose $g(t) \in \mathcal{G}^\C$ such that $g(t) g^*(t) = h(t)^{-1}$ (Note that \emph{a priori} this choice is not unique). 

\begin{lem}\label{lem:equivflowtogaugeflow}
Let $h(t)$ be a solution to equation \eqref{eqn:gaugeflow}, choose $g(t) \in \mathcal{G}^\C$ such that $g(t) g(t)^* = h(t)^{-1}$, and let $A''(t) = g(t) \cdot A_0''$, $\phi(t) = g(t) \cdot \phi_0$. Then  $(A''(t), \phi(t))$ is a solution to \eqref{eqn:equivflow} with $\alpha(t) = \frac{1}{2} (g^{-1} \partial_t g - (\partial_t g^*)(g^*)^{-1} )$.
\end{lem}

\begin{proof}[Proof of Lemma \ref{lem:equivflowtogaugeflow}]
Let $(A''(t), \phi(t)) = (g(t) \cdot A_0'', g(t) \cdot \phi_0)$. We have the following identities for $g \in \mathcal{G}^\C$ (cf \cite{Hong01} $(3.2)$ for the vortex equations)
\begin{equation}\label{eqn:curvaturetransform}
g F_A g^{-1} = g F_{(g \cdot A_0)} g^{-1} = F_{A_0} + d_{A_0}''(h^{-1} (d_{A_0}' h))
\end{equation}
\begin{equation}\label{eqn:phibrackettransform}
g [\phi, \phi^*] g^{-1} = g [(g\cdot \phi_0), (g \cdot \phi_0^*)] g^{-1} = [\phi_0, h^{-1} \phi_0^* h]
\end{equation}
Differentiating $A''$ and $\phi$ gives us
\begin{align}\label{eqn:Aderivative}
\begin{split}
\frac{\partial A''}{\partial t} & = \left. \frac{\partial}{\partial \varepsilon} \right|_{\varepsilon=0} d_{(g+\varepsilon \partial_t g)\cdot A_0}'' = d_A''(g^{-1}(\partial_t g)) \\
 & = \frac{1}{2} d_A''({g}^{-1} \partial_t g + (\partial_t g^*)(g^*)^{-1} ) + \frac{1}{2} d_A'' (g^{-1} \partial_t g - (\partial_t g^*)(g^*)^{-1} )
\end{split}
\end{align}
and similarly
\begin{equation}\label{eqn:phiderivative}
\frac{\partial \phi}{\partial t} =  \frac{1}{2} [\phi, (g^{-1} \partial_t g + (\partial_t g^*)(g^*)^{-1} )] + \frac{1}{2} [\phi, (g^{-1} \partial_t g - (\partial_t g^*)(g^*)^{-1} )]
\end{equation}
Let $\alpha(t) = \frac{1}{2} (g^{-1} \partial_t g - (\partial_t g^*)(g^*)^{-1} )$. Since $g g^* = h^{-1}$, then
\begin{equation}\label{eqn:hderivative}
\frac{\partial h}{\partial t} = -(g^*)^{-1}\left((\partial_t g^*) (g^*)^{-1} + g^{-1} (\partial_t g) \right)g^{-1}
\end{equation}
Using the identities \eqref{eqn:curvaturetransform} and \eqref{eqn:phibrackettransform} together with the equation \eqref{eqn:gaugeflow} shows that the right-hand side of \eqref{eqn:hderivative} is $-2i (g^*)^{-1} g^{-1} g *(F_A + [\phi, \phi^*] ) g^{-1} + 2i \lambda h$, and therefore
\begin{equation*}
\frac{1}{2}\left((\partial_t g^*) (g^*)^{-1} + g^{-1} (\partial_t g)\right) = i*(F_A + [\phi, \phi^*]) - i \lambda \cdot \id
\end{equation*}
Together with \eqref{eqn:Aderivative} and \eqref{eqn:phiderivative} this gives us the following equations for $A''(t)$ and $\phi(t)$
\begin{align*}
\frac{\partial A''}{\partial t} & = id_A'' * (F_A + [\phi, \phi^*]) + d_A'' (\alpha - i\lambda \cdot \id) \\
\frac{\partial \phi}{\partial t} & = i[\phi, (F_A + [\phi, \phi^*])] + [\phi, \alpha - i\lambda \cdot \id] \qedhere
\end{align*}

\end{proof}

\begin{proof}[Proof of Proposition \ref{lem:existence}]

To prove existence, we construct a solution to the gradient flow equations \eqref{eqn:gradfloweqns} from a solution to the equivalent flow equations \eqref{eqn:equivflow}. To prove uniqueness we then show that this solution is independent of the choice of $g(t)$ such that $g(t) g(t)^* = h(t)^{-1}$. Consider the following ODE for a one-parameter family of complex gauge transformations $S(t)$
\begin{equation}\label{eqn:equivgaugetrans}
\frac{\partial S}{\partial t} = S(t) \left( \alpha(t) - i \lambda \cdot \id\right) 
\end{equation}
where $\alpha : \R \rightarrow \Lie(\mathcal{G})$ is as defined in the proof to Lemma \ref{lem:equivflowtogaugeflow}. Note firstly that $S(t)$ is a unitary gauge transformation, even though \emph{a priori} $S(t) \in \mathcal{G}^\C$. This follows from observing that $S(0)=\id \in \mathcal{G}$ and $\frac{\partial S}{\partial t} \in S(t) \cdot \Lie(\mathcal{G})$, therefore $S(t) \in \mathcal{G}$ for all $t$. Lemma \ref{lem:equivflowtogaugeflow} shows that $\alpha(t)$ is defined for all $t$, and therefore solutions to equation \eqref{eqn:equivgaugetrans} exist for all time by linear ODE theory. 

Let $\left( \tilde{A}(t), \tilde{\phi}(t) \right)$ denote a solution to the equivalent flow equations. For notation let $\tilde{\alpha} = \alpha - i\lambda \cdot \id$. Define $A''(t) = S^{-1}(t) \cdot \tilde{A''}(t)$ and $\phi(t) = S(t)^{-1} \cdot \tilde{\phi}(t)$. Then $\left(A''(t), \phi(t) \right)$ exists for all $t$ and it remains to show that $\left(A''(t), \phi(t) \right)$ satisfies the gradient flow equations \eqref{eqn:gradfloweqns}. Differentiating with respect to $t$ gives us
\begin{align*}
\frac{\partial A''}{\partial t}  & = \frac{\partial}{\partial t} \left( S \tilde{d_A''}  S^{-1}  \right) \\
 & = S\tilde{\alpha} \tilde{d_A''} S^{-1}  + S i\tilde{d_A''} * \left( \tilde{F_A} + [\tilde{ \phi}, \tilde{ \phi^*} ] \right) S^{-1} + S \left( \tilde{d_A''} \tilde{\alpha} \right) S^{-1}  - d_A''  \frac{\partial S}{\partial t} S^{-1}  \\
 & = S \tilde{\alpha } \tilde{d_A''} S^{-1} + i d_A'' * \left( F_A + [\phi, \phi^*] \right) + S \left( \tilde{d_A''} \tilde{\alpha} \right) S^{-1}  -S \tilde{d_A''} \tilde{\alpha} S^{-1} \\
 & = i d_A'' * \left( F_A + [\phi, \phi^*] \right)
\end{align*}
and similarly for $\tilde{\phi}$ we obtain $\frac{\partial \phi}{dt} = i[\phi, *(F_A + [\phi, \phi^*])]$. Therefore the solution $\left( A''(t), \phi(t) \right)$ of \eqref{eqn:gradfloweqns} exists for all time.

To prove uniqueness we note firstly that (as in the unitary case studied in \cite{Donaldson85}) solutions to Simpson's heat equation \eqref{eqn:simpsonheat} are unique, by applying the maximum principle to the distance function $\sigma$ given in the proof of Proposition 6.3 of \cite{Simpson88}. From the construction in the proofs of Lemmas \ref{lem:gaugeflowtosimpson} and \ref{lem:equivflowtogaugeflow}, the only non-unique choice made in constructing the solution to the gradient flow of $\YMH$ from a solution to equation \eqref{eqn:simpsonheat} is the choice of $g(t)$ such that $g(t) g(t)^* = h(t)^{-1}$. The following lemma shows that the solution is independent of this choice.

\begin{lem}\label{lem:uniquenessofgS}
Let $h(t)$ be a solution to \eqref{eqn:simpsonheat}, and suppose that $g_1(t)$ and $g_2(t)$ are one parameter families in $\mathcal{G}^\C$ such that $g_1(t) g_1(t)^* = h(t)^{-1} = g_2(t) g_2(t)^*$. Let $S_1(t)$ and $S_2(t)$ be the corresponding solutions constructed above such that
\begin{align*}
\left( A_1''(t), \phi_1(t) \right) & = S_1(t)^{-1} \cdot g_1(t) \cdot \left( A_0'', \phi_0 \right) = \left( g_1(t) S_1(t)^{-1} \right) \cdot  \left( A_0'', \phi_0 \right) \\
\left( A_2''(t), \phi_2(t) \right) & = S_2(t)^{-1} \cdot g_2(t) \cdot \left( A_0'', \phi_0 \right) = \left( g_2(t) S_2(t)^{-1} \right) \cdot  \left( A_0'', \phi_0 \right) 
\end{align*}

Then $\left( A_1''(t), \phi_1(t) \right) = \left( A_2''(t), \phi_2(t) \right)$.

\end{lem}

\begin{proof}[Proof of Lemma \ref{lem:uniquenessofgS}]

Note that $g_1^{-1} g_2 g_2^* {(g_1^*)}^{-1} = \id$, therefore $g_1^{-1} g_2 = u(t)$ for some curve $u(t) \in \mathcal{G}$. As in the proof of Lemma \ref{lem:equivflowtogaugeflow}, define the gauge fixing terms $\alpha_1(t)$ and $\alpha_2(t)$ by
\begin{align*}
\alpha_1(t) & = \frac{1}{2} \left( g_1^{-1} \partial_t g_1 - ( \partial_t g_1^* ) {(g_1^*)}^{-1} \right) \\
\alpha_2(t) & = \frac{1}{2} \left( g_2^{-1} \partial_t g_2 - ( \partial_t g_2^* ) {(g_2^*)}^{-1} \right) = u(t)^{-1} \alpha_1(t) u(t) + u(t)^{-1} \partial_t u
\end{align*}

Therefore the equations for $S_1(t)$ and $S_2(t)$ are
\begin{align*}
S_1(t)^{-1} \frac{\partial S}{\partial t} & = \alpha_1(t) - i\lambda \\
S_2(t)^{-1} \frac{\partial S}{\partial t} & = \alpha_2(t) - i\lambda = u(t)^{-1} \alpha_1(t) u(t) + u(t)^{-1} \partial_t u
\end{align*}

$S_2(t) = S_1(t) u(t)$  is a solution to this equation, which is unique by linear ODE theory. Therefore $g_2(t) S_2(t)^{-1} = g_1(t) u(t) u(t)^{-1} S_1(t)^{-1} = g_1(t) S_1(t)^{-1}$, which completes the proof of uniqueness.
\end{proof} \renewcommand{\qedsymbol}{} \end{proof}

\subsection{Compactness along the gradient flow}\label{subsec:compactness}

In this section we derive estimates for $\left| \nabla_A^k (F_A + [\phi, \phi^*]) \right|_{C^0}$ along the gradient flow of $\YMH$, and prove a compactness theorem  (Lemma \ref{lem:L4compactness}). Together these are sufficient to prove Proposition \ref{lem:convmodgauge}. The basic tool is the following estimate based on Theorem 4.3 in \cite{Hitchin87} (for the case of $\SU(2)$ bundles) and Lemma 2.8 of \cite{Simpson92} (for bundles with a general compact structure group).

\begin{thm}[Hitchin/Simpson]\label{thm:hitchinuhlen}
Fix a Higgs pair $(A_0, \phi_0)$ and a constant $C$, and consider the subset $\mathcal{O}_C$ of the complex group orbit $\mathcal{G}^\C \cdot (A_0, \phi_0)$ consisting of Higgs pairs satisfying the estimate $\| F_A + [\phi, \phi^*] \|_{L^2} < C$. Then there exists a constant $K$ such that $\| F_A \|_{L^2} < K$ and $\| \phi \|_{H^1} < K$ for all $(A, \phi) \in \mathcal{O}_C$.
\end{thm}

The Sobolev spaces $L_k^p$ used in this section are defined via norm
\begin{equation*}
\| \sigma \|_{L_k^p}^A = \left( \sum_{i=0}^k \| \nabla_A^i \sigma \|_{L^p} \right).
\end{equation*}

\begin{rek}
\emph{A priori} the norm depends on the connection $d_A$, however Propositions D.1 and D.2 from \cite{radethesis} show that given a uniform bound on the curvature $\| F_A \|_{L_k^2}^A$, the norms of the Sobolev multiplication, embedding and interpolation operators are uniformly bounded in $A$. Therefore the bounds obtained from Lemma \ref{lem:L4compactness} below show that the estimates obtained in this section are independent of the choice of connection used to define the Sobolev norm.
\end{rek}

The proof of Proposition \ref{lem:convmodgauge} relies on the following two lemmas. Firstly, by bootstrapping the results of Theorem \ref{thm:hitchinuhlen} using the equation $d_A'' \phi = 0$ we obtain the following result.

\begin{lem}\label{lem:L4compactness}
Consider the subset $\mathcal{O}_C^k$ of the complex group orbit $\mathcal{G}^\C \cdot (A_0, \phi_0)$ consisting of Higgs pairs satisfying the estimate $\| F_A + [\phi, \phi^*] \|_{L_k^4} < C$. Then there exists a constant $K$ such that $\| F_A \|_{L_k^4} < K$ and $\| \phi \|_{L_{k+2}^4} < K$ for all $(A'', \phi) \in \mathcal{O}_C^k$. Moreover, the Sobolev embedding theorems show that $\| \nabla_A^{k+1} \phi \|_{C^0} < K$.
\end{lem}

\begin{proof}[Proof of Lemma \ref{lem:L4compactness}]
Suppose that $\| F_A + [\phi, \phi^*] \|_{L_k^4} < C$ on a $\mathcal{G}^\C$-orbit. Then $\| F_A + [\phi, \phi^*] \|_{L^2}$ is bounded and Theorem \ref{thm:hitchinuhlen} shows that there exists $K$ such that $\| F_A \|_{L^2} < K$ and $\| \phi \|_{L_1^2} < K$. Therefore $[\phi, \phi^*]$ is bounded in $L^4$ and so $\| F_A \|_{L^4}$ is bounded. Theorem 1.5 in \cite{Uhlenbeck82} shows that after applying unitary gauge transformations $\| A \|_{L_1^4} < K$ locally. Then Sobolev multiplication $L_1^4 \times L_1^2 \rightarrow L_1^2$ shows that $\| [A'', \phi] \|_{L_1^2}$ is bounded locally and so the equation $d'' \phi = - [A'', \phi]$ gives the elliptic estimate $\| \phi \|_{L_2^2} \leq C(\| [A'', \phi] \|_{L_1^2} + \| \phi \|_{L^2})$. Applying this procedure again with Sobolev multiplication $L_1^4 \times L_2^2 \rightarrow L_1^4$ shows that $\phi$ is bounded in $L_2^4$. Therefore $[\phi, \phi^*]$ is bounded in $L_1^4$, so $\| F_A \|_{L_1^4} < C$ and we can repeat the above process inductively for all $k$ to complete the proof of Lemma \ref{lem:L4compactness}.
\end{proof}

The next lemma shows that the $L_k^4$ bound on $F_A + [\phi, \phi^*]$ exists along the flow.
\begin{lem}\label{lem:chenshen}
Let $s \geq 0$ and suppose that $\| \nabla_A^\ell (F_A + [\phi, \phi^*]) \|_{C^0}$ is bounded for all $\ell < s$, and that $\| \nabla_A^\ell \phi \|_{C^0}$ is bounded for all $\ell \leq s$. Then the following estimates hold for a solution $(A(t), \phi(t))$ of the gradient flow equations \eqref{eqn:gradfloweqns}
\begin{multline}\label{eqn:chensheninductivebound}
2 \left| \nabla_A^s (F_A + [\phi, \phi^*]) \right|^2  \leq - \left( \frac{\partial}{\partial t} + \Delta \right) \left| \nabla_A^{s-1} (F_A + [\phi, \phi^*] ) \right|^2 \\
  + C \left( \left| \nabla_A^{s-1} (F_A + [\phi, \phi^*]) \right|^2 + 1 \right) 
\end{multline}
\begin{equation}\label{eqn:chenshenparaboliceqn}
\left( \frac{\partial}{\partial t} + \Delta \right) \left| \nabla_A^s (F_A + [\phi, \phi^*]) \right|^2  \leq C \left( \left| \nabla_A^s (F_A + [\phi, \phi^*]) \right|^2 + 1 \right) 
\end{equation}
\end{lem}

\begin{proof}[Proof of Lemma \ref{lem:chenshen}]
For notation let $\mu = F_A + [\phi, \phi^*]$ and define the operator $L : \Omega^0(\ad(E)) \rightarrow \Omega^{1, 0}(\End(E))$ by $L(u) = [\phi, u]$. Firstly we note that (in the notation of Section \ref{sec:symplectic}) for any moment map $\mu$ on a symplectic manifold we have the following equation along the gradient flow
\begin{equation}
\frac{\partial (*\mu)}{\partial t} = *d\mu \left(\frac{\partial x}{\partial t} \right) = -\rho_x^* I (-I \rho_x(*\mu)) = -\rho_x^* \rho_x (*\mu)
\end{equation}
For Higgs bundles this reduces to the equation $\left( \frac{\partial}{\partial t} + \Delta_{(A'', \phi)} \right)(*\mu) = 0$. Since $*\mu$ is a $0$-form then $\Delta_{(A'', \phi)}(*\mu) = \nabla_A^* \nabla_A (*\mu) + L^* L (*\mu)$. The method of \cite{Donaldson85} pp16-17 for the Yang-Mills functional shows that in this case
\begin{equation}\label{eqn:normcurvatureheat}
\frac{\partial \left| *\mu \right|^2}{\partial t} + \Delta \left| * \mu \right|^2 \leq 0
\end{equation}
In particular, the maximum principle shows that $\sup_X \left| * \mu \right|^2$ is decreasing and therefore bounded uniformly in $t$. Equations \eqref{eqn:chensheninductivebound} and \eqref{eqn:chenshenparaboliceqn} can then be computed in a standard way (cf p40 of \cite{ChenShen93} for the Yang-Mills flow and the proof of Proposition 3 and Proposition 6 in \cite{HongTian04} for the vortex equations), and so the rest of the proof is omitted. 
\end{proof}

As a corollary, we obtain uniform $L_k^2$ bounds on $F_A + [\phi, \phi^*]$. 
\begin{cor}\label{cor:curvaturebounds}
$\displaystyle{\int_T^{T+1} \| \nabla_A^s (F_A + [\phi, \phi^*]) \|_{L^2} \, dt}$ is bounded uniformly in $T$, and so $\| \nabla_A^s (F_A + [\phi, \phi^*]) \|_{C^0}$ is bounded uniformly in $t$.
\end{cor}

The proof relies on Moser's Harnack inequality from \cite{Moser64}, which can be stated in the following form.
\begin{thm}[Moser]\label{thm:moserharnack}
Let $0 \leq \tau_1^- < \tau_2^- < \tau_1^+ < \tau_2^+$ and suppose that $u \geq 0$ is a function on a compact manifold $X$, and that $\frac{\partial u}{\partial t} + \Delta u \leq C u$. Then there exists a constant $\gamma$ depending only on $(\tau_2^- - \tau_1^-)$, $(\tau_1^+ - \tau_2^-)$, $(\tau_2^+ - \tau_2^-)$ and $C$ such that
\begin{equation*}
\sup_{\tau_1^- < t < \tau_2^-} u \leq \gamma \int_{\tau_1^+}^{\tau_2^+} \| u \|_{L^2} \, dt
\end{equation*}
\end{thm}

\begin{proof}[Proof of Corollary \ref{cor:curvaturebounds}]
 To obtain a $C^0$ bound on $\left| \nabla_A^s (*\mu) \right|$ we use Theorem \ref{thm:moserharnack} as follows. Equation \eqref{eqn:chensheninductivebound} together with the fact that $\left| \nabla_A^\ell (*\mu) \right|$ is bounded in $C^0$ for all $\ell < s$ shows that $\int_T^{T+1} \| \nabla_A^s (*\mu) \|_{L^2} \, dt < C$, where $C$ is independent of $T$. Equation \eqref{eqn:chenshenparaboliceqn} shows that Moser's theorem applies to the function $\left| \nabla_A^s (*\mu) \right| + 1$. Therefore
\begin{equation*}
\sup_{T-2 < t < T-1} \left| \nabla_A^s (*\mu) \right| + 1 \leq \gamma \int_T^{T+1} \| \nabla_A^s (*\mu) \|_{L^2} \, dt
\end{equation*}
is uniformly bounded in $T$ (where $\gamma$ is independent of $T$ because the time intervals $[T-2, T-1]$ and $[T, T+1]$ are of constant size and relative position). Therefore $\left| \nabla_A^s (*\mu) \right|$ is uniformly bounded in $t$.
\end{proof}

Using these lemmas, the proof of Proposition \ref{lem:convmodgauge} proceeds as follows.

\begin{proof}[Proof of Proposition \ref{lem:convmodgauge}]

Firstly we show by induction that $\| F_A + [\phi, \phi^*]  \|_{L_k^4}$ is bounded for all $k$. The computation in the proof of Lemma \ref{lem:chenshen} shows that a solution $(A(t), \phi(t))$ of the gradient flow equations \eqref{eqn:gradfloweqns} satisfies the equation $$\left( \frac{\partial}{\partial t} + \Delta \right)\left| F_A + [\phi, \phi^*] \right|^2 \leq 0$$ Therefore $\| F_A + [\phi, \phi^*] \|_{C^0}$ is bounded uniformly in $t$, and in particular $\| F_A + [\phi, \phi^*] \|_{L^4}$ is bounded. Lemma \ref{lem:L4compactness} then gives a bound on the $C^0$ norm of $\left| \nabla_A \phi \right|$ and Corollary \ref{cor:curvaturebounds} with $s=1$ gives a bound on $\| \nabla_A(F_A + [\phi, \phi^*]) \|_{C^0}$, and hence on $\| F_A + [\phi, \phi^*] \|_{L_1^4}$.

Now suppose that $\| F_A + [\phi, \phi^*] \|_{L_k^4}$ is bounded, and also suppose that $\| \nabla_A^\ell \phi \|_{C^0}$ and $\| \nabla_A^\ell \left( F_A + [\phi, \phi^*] \right) \|_{C^0}$ are bounded uniformly in $t$ for all $\ell \leq k$. Applying Lemma \ref{lem:L4compactness} shows that $\| \nabla_A^\ell \phi \|_{C^0}$ is bounded for all $\ell \leq k+1$. Then we can apply Lemma \ref{lem:chenshen} for $s = k+1$ which shows that $\| F_A + [\phi, \phi^*] \|_{L_{k+1}^4}$ and $\| \nabla_A^{k+1} (F_A + [\phi, \phi^*]) \|_{C^0}$ are bounded, which completes the induction.

Since $\| (F_A + [\phi, \phi^*]) \|_{L_k^4}$ is bounded for all $k$ then Lemma \ref{lem:L4compactness} holds for all $k$. In particular, $\| F_A \|_{L_k^4}$ and $\| \phi \|_{L_{k+2}^4}$ are bounded for all $k$. To complete the proof we need to show that along a subsequence the gradient flow converges to a critical point of $\YMH$. To see this, firstly note that in general for the gradient flow of any non-negative functional $f : M \rightarrow \R$  we have for any time $T$ the equation $f(t=0) - f(t=T) = -\int_0^T \frac{\partial f}{\partial t} \, dt = \int_0^T df (\grad f) \, dt$ and therefore $\int_0^T \| \grad f \|^2 \, dt \leq f(t=0)$. Therefore there exists a subsequence $t_n \rightarrow \infty$ such that $\grad f(t_n) \rightarrow 0$ strongly in the appropriate norm. For the case of $f = \YMH$, along this subsequence $t_n$ the above argument provides a bound on $\| F_A \|_{L_k^4}$. Therefore Uhlenbeck's compactness theorem shows that along a subsequence (also call it $t_n$) there exists a sequence of unitary gauge transformations $g_n$ such that $g_n \cdot A(t_n) \rightharpoonup A_\infty$ weakly in $L_{k+1}^4$ and strongly in $L_k^4$. Since $\| g_n \cdot \phi(t_n) \|_{L_{k+2}^4}$ is also bounded, then there exists a subsequence (also call it $t_n$) such that $g_n \cdot \phi(t_n) \rightarrow \phi_\infty$ in $L_k^4$. It only remains to show that $(A_\infty, \phi_\infty)$ is a critical point of $\YMH$. 

Let $\rho_n: \Omega^0(\ad(E)) \rightarrow \Omega^{0,1}(\End(E)) \oplus \Omega^{1,0}(\End(E))$ denote the operator 
$$u \mapsto \left( \begin{matrix} d_{A(t_n)}'' u \\ [\phi(t_n), u] \end{matrix} \right)$$
and let $*\mu = *(F_A + [\phi, \phi^*])$. Note that $\grad \YMH(t_n) = I \rho_n(*\mu(t_n))$. Along the subsequence $t_n$, $\grad \YMH \rightarrow 0$ strongly in $L_{k-1}^4$. Therefore 
\begin{equation}
\rho_n(*\mu(t_n)) - \rho_\infty(*\mu(\infty)) = \rho_\infty(*\mu(t_n) - *\mu(\infty)) - (\rho_n - \rho_\infty)(*\mu(t_n))
\end{equation}
$\rho_n(*\mu(t_n)) \rightarrow 0$ strongly in $L_{k-1}^4$ and the right-hand side of the above equation converges to $0$ strongly in $L_{k-1}^4$. Therefore $\rho_\infty(*\mu(\infty)) = 0$, and so $(A_\infty, \phi_\infty)$ is a critical point of $\YMH$. 
\end{proof}

\subsection{Continuous Dependence on Initial Conditions}

In this section we prove Proposition \ref{lem:continuousdependence}. The proof of this proposition follows the method of Section 5 in \cite{Rade92} which proves continuous dependence on the initial conditions in the $H^1$ norm for the Yang-Mills gradient flow, however here we generalise to the case of Higgs bundles, and also use the estimates for the higher derivatives of the curvature from Lemma \ref{lem:L4compactness} to show continuous dependence on the initial conditions in the $H^k$ norm for all $k$. This relies on the estimates from Proposition A of \cite{Rade92}, which are valid when the higher derivatives of the curvature are bounded. Rade's approach also proves the existence and uniqueness of a solution, however since in this case Proposition \ref{lem:existence} together with the estimates derived in the proof of Proposition \ref{lem:convmodgauge} already show that a unique smooth solution to \eqref{eqn:gradfloweqns} exists, then the estimates in this section can be simplified from those of Section 5 in \cite{Rade92}. The reference for the definitions of the time-dependent Sobolev spaces used in this section is the Appendix of \cite{Rade92} (further details are explained in \cite{radethesis}). 

Firstly note that for the general case of a moment map on a symplectic manifold, the downwards gradient flow of $\| \mu \|^2$ satisfies the following equations
\begin{align}\label{eqn:generalmomentmapgradfloweqns}
\begin{split}
\frac{\partial x}{\partial t} + I \rho_x(*\mu) & = 0 \\
\frac{\partial (*\mu)}{\partial t} + \rho_x^* \rho_x (*\mu) & = 0
\end{split}
\end{align}
The results of Proposition \ref{lem:existence} and Proposition \ref{lem:convmodgauge} show that in the Higgs bundle case, a smooth solution to \eqref{eqn:generalmomentmapgradfloweqns} exists. Now consider instead the following generalised system, with $*\mu$ replaced by a general $\Omega \in \Lie(\mathcal{G})$
\begin{align}\label{eqn:generalgradfloweqns}
\begin{split}
\frac{\partial x}{\partial t} + I \rho_x(\Omega) & = 0 \\
\frac{\partial \Omega}{\partial t} + \rho_x^* \rho_x (\Omega) & = 0
\end{split}
\end{align}

Firstly we note that if a smooth solution $(x(t), \Omega(t))$ of \eqref{eqn:generalgradfloweqns} exists with initial conditions $x(0) = x_0$ and $\Omega(0) = *\mu(x_0)$ then this solution satisfies $\Omega(t) = *\mu(x(t))$. This follows by considering $\psi(t) = \Omega(t) - *\mu(x(t))$, and noting that
\begin{align*}
\frac{\partial \psi}{\partial t} & = \frac{\partial \Omega}{\partial t} - \frac{\partial (*\mu) }{\partial t} \\
 & = -\rho_x^* \rho_x(\Omega) - * d\mu \left( \frac{\partial x}{\partial t} \right) \\ 
 & = -\rho_x^* \rho_x(\Omega) + \rho_x^* I \left( -I \rho_x(\Omega) \right) \\
 & = -\rho_x^* \rho_x(\Omega) + \rho_x^* \rho_x (\Omega) = 0
\end{align*}
Therefore if $\psi(t=0) = 0$ then $\Omega(t) = *\mu(x(t))$ for all $t$. In the Higgs bundle case, the space $T^* \mathcal{A}$ is an affine space, and $\rho_{x+a}(u) = \rho_x(u) + \{ a, u \}$ where $\{ \cdot, \cdot \}$ denotes various intrinsically defined multilinear operators. For a fixed point $x_0 \in T^* \mathcal{A}$, let $y = x - x_0$ and note that the equations \eqref{eqn:generalgradfloweqns} become
\begin{align}\label{eqn:expandedgradfloweqns}
\begin{split}
\frac{\partial y}{\partial t} + I \rho_{x_0}(\Omega) & = \{ y , \Omega \} \\
\frac{\partial \Omega}{\partial t} + \rho_{x_0}^* \rho_{x_0}( \Omega) & = \{ y^*, \rho_{x_0}(\Omega) \} + \{ \rho_{x_0}^* (y), \Omega \} + \{ y^*, y, \Omega \}
\end{split}
\end{align}
In the Higgs bundle case we can write (for $x_0 = \left( A_0, \phi_0 \right)$)
\begin{align*}
\rho_{x_0}^* \rho_{x_0} (\Omega) & = d_{A_0}''^* d_{A_0}'' \Omega - \starbar \left[\phi_0, \starbar[\phi_0, \Omega] \right] \\
 & = \nabla_{A_0}^* \nabla_{A_0} \Omega + \{ \phi_0^*, \phi_0, \Omega \}
\end{align*}
Therefore the gradient flow equations become
\begin{align}
\begin{split}
\frac{\partial y}{\partial t} + I \rho_{x_0}(\Omega) & = \{ y, \Omega \} \\
\frac{\partial \Omega}{\partial t} + \nabla_{A_0}^* \nabla_{A_0} \Omega & = \{ \phi_0^*, \phi_0, \Omega \} + \{y^*, \rho_{x_0}(\Omega) \} + \{ y^*, y, \Omega \}
\end{split}
\end{align}
Following the method of \cite{Rade92}, define the operator $L$
\begin{equation}
L = \left( \begin{matrix} \frac{\partial }{\partial t} & \rho_{x_0}^* \\  0 & \frac{\partial}{\partial t} + \nabla_{A_0}^* \nabla_{A_0} \end{matrix} \right)
\end{equation}
and $Q_1$, $Q_2$, $Q_3$
\begin{align*}
Q_1 \left( \begin{matrix} y \\ \Omega \end{matrix} \right) & = \left( \begin{matrix} 0 \\ \{ \phi_0^*, \phi_0, \Omega \} \end{matrix} \right) \\
Q_2 \left( \begin{matrix} y \\ \Omega \end{matrix} \right) & = \left( \begin{matrix} \{ y, \Omega \} \\ \{ y, \rho_{x_0}(\Omega) \} + \{ \rho_{x_0}^* (y), \Omega \} \end{matrix} \right) \\
Q_3 \left( \begin{matrix} y \\ \Omega \end{matrix} \right) & = \left( \begin{matrix} 0 \\ \{ y^*, y, \Omega \} \end{matrix} \right) 
\end{align*}

Define the Hilbert spaces
\begin{align*}
\mathcal{U}^k(t_0) & = \left\{ \left( \begin{matrix} y \\ \Omega \end{matrix} \right) : y \in H^{\frac{1}{2}+\varepsilon, k}([0, t_0]) \, \mathrm{ and } \, \Omega \in H^{\frac{1}{2}+\varepsilon, k-1}([0, t_0]) \cap H^{-\frac{1}{2} + \varepsilon, k+1}([0, t_0]) \right\} \\
\mathcal{U}_P^k(t_0) & = \left\{ \left( \begin{matrix} y \\ \Omega \end{matrix} \right) : y \in H_P^{\frac{1}{2}+\varepsilon, k}([0, t_0]) \, \mathrm{ and } \, \Omega \in H_P^{\frac{1}{2}+\varepsilon, k-1}([0, t_0]) \cap H_P^{-\frac{1}{2} + \varepsilon, k+1}([0, t_0]) \right\} \\
\mathcal{W}_P^k(t_0) & = \left\{ \left( \begin{matrix} y \\ \Omega \end{matrix} \right) : y \in H_P^{-\frac{1}{2}+\varepsilon, k}([0, t_0]) \, \mathrm{ and } \, \Omega \in H_P^{-\frac{1}{2}+\varepsilon, k-1}([0, t_0])  \right\}
\end{align*}

The following lemma is a Higgs bundle-version of \cite{Rade92} Lemma 5.1, the proof is analogous and therefore omitted.
\begin{lem}\label{lem:bounds}
Let $(A_0, \phi_0) \in \mathcal{B}$. Then the maps $L$, $Q_i$ for $i = 1,2,3$ and the identity map $I$ define bounded linear operators
\begin{align*}
L : & \mathcal{U}_P^k(t_0) \rightarrow \mathcal{W}_P^k(t_0) \\
Q_1 : & \mathcal{U}^k(t_0) \rightarrow \mathcal{W}_P^k(t_0) \\
Q_2 : & S^2 \mathcal{U}^k(t_0) \rightarrow \mathcal{W}_P^k(t_0) \\
Q_3 : & S^3 \mathcal{U}^k(t_0) \rightarrow \mathcal{W}_P^k(t_0) \\
I : & \mathcal{U}_P^k(t_0) \rightarrow \mathcal{U}(t_0)
\end{align*}
Moreover, the operator $L$ is invertible. For any $K>0$ there exists $c_K > 0$ such that if $\| F_{A_0} \|_{H^{k-1}} < K$ then 
\begin{align*}
\| L^{-1} \| \leq c_K, \quad & \| Q_1 \| \leq c_K t_0^{\frac{1}{4} - \varepsilon} \\
\| M \| \leq c_K t_0^{-\varepsilon}, \quad & \| Q_2 \| \leq c_K t_0^{\frac{1}{4} - 2 \varepsilon} \\
\| I \| \leq 1, \quad & \| Q_3 \| \leq c_K t_0^{\frac{1}{2} - 2 \varepsilon}
\end{align*}
for $t_0$ sufficiently small
\end{lem}
Note that the Sobolev spaces in \cite{Rade92} are defined slightly differently to the definitions of $\mathcal{U}^k$, $\mathcal{U}_P^k$ and $\mathcal{W}^k$ above. Rade also considers the case of a three-dimensional manifold for which the multiplication theorems used in the proof of Lemma 5.1 of \cite{Rade92} become borderline with the definitions above, however here we only consider the case of a compact Riemann surface, and so we can derive stronger estimates. Now consider the homogeneous system of equations with initial conditions $(y_1(0), \Omega_1(0)) = (x_0, \Omega_0)$.
\begin{align}
\frac{\partial y_1}{\partial t} + I \rho_{x_0}(\Omega_1) & = 0 \label{eqn:homogeneousy} \\
\frac{\partial \Omega_1}{\partial t} + \nabla_{A_0}^* \nabla_{A_0} \Omega_1 & = 0 \label{eqn:homogeneousomega}
\end{align}
Proposition A of \cite{Rade92} shows that there exists a unique solution to \eqref{eqn:homogeneousomega} given by $\Omega_1 \in H^{\frac{1}{2} + \varepsilon, k-1} \cap H^{-\frac{1}{2} + \varepsilon, k+1}$, which satisfies $\| \Omega_1 \| \leq c_K t_0^{-\varepsilon} \| \Omega_0 \|$. Therefore there also exists a unique solution $y_1 = y_0 - \int_0^t I \rho_{x_0}(\Omega_1(s)) \, ds$ to \eqref{eqn:homogeneousy} which (again by Proposition A of \cite{Rade92}) satisfies $\| y_1 \| \leq c_K \| y_0 \|$. Therefore the solution operator $M$ defined by
\begin{equation}
M\left( \begin{matrix} y_0 \\ \Omega_0 \end{matrix} \right) = \left( \begin{matrix} y_1 \\ \Omega_1 \end{matrix} \right)
\end{equation}
is bounded, with $\| M \| \leq c_K t_0^{-\varepsilon}$. Let $(y_2, \Omega_2) = (y-y_1, \Omega - \Omega_1)$. Then the initial-value problem \eqref{eqn:expandedgradfloweqns} can be written as
\begin{multline}
L\left( \begin{matrix} y_2 \\ \Omega_2 \end{matrix} \right) = Q_1 \left( M\left( \begin{matrix} y_0 \\ \Omega_0 \end{matrix} \right) + I \left( \begin{matrix} y_2 \\ \Omega_2 \end{matrix} \right) \right) + Q_2 \left( M \left( \begin{matrix} y_0 \\ \Omega_0 \end{matrix} \right) + I \left( \begin{matrix} y_2 \\ \Omega_2 \end{matrix} \right) \right) \\
 + Q_3 \left( M \left( \begin{matrix} y_0 \\ \Omega_0 \end{matrix} \right) + I \left( \begin{matrix} y_2 \\ \Omega_2 \end{matrix} \right) \right)
\end{multline}

The estimates from Lemma \ref{lem:bounds} are identical to those of Lemma 5.1 in \cite{Rade92}, and applying Lemma 5.2 of \cite{Rade92} shows that for a small interval $[0, t_0]$, the solution to \eqref{eqn:expandedgradfloweqns} satisfies $y \in C^0([0, t_0], H^k)$, $\Omega \in C^0([0, t_0], H^{k-1})$, and that $(y, \Omega)$ depends continuously on the initial conditions $(y_0, \Omega_0) \in H^k \times H^{k-1}$. This completes the proof of Proposition \ref{lem:continuousdependence}.

\subsection{A Lojasiewicz inequality}\label{subsec:loja}

In the paper \cite{Simon83}, Simon proved the convergence of solutions to the equation
$$
\dot{u} - \mathcal{M}(u) = f
$$
as $t \rightarrow \infty$, where $u = u(x, t)$ is a smooth section of a vector bundle $F$ over a compact Riemannian manifold $\Sigma$, and $\mathcal{M}(u)$ is the gradient of an "Energy Functional" $\mathcal{E}(u) = \int_{\Sigma} E(x, u, \nabla u)$ on $\Sigma$. The function $E$ is assumed to have analytic dependence on $u$ and $\nabla u$, and the operator $\mathcal{M}$ is assumed to be elliptic. The key estimate in Simon's proof was the inequality
$$
\left\| \mathcal{M}(u) \right\| \geq \left| \mathcal{E}(u) - \mathcal{E}(0) \right|^{1-\theta}
$$
where $\theta \in (0, \frac{1}{2})$, an infinite dimensional version of an inequality proved by Lojasiewicz in \cite{Lojasiewicz84} for real analytic functionals on a finite-dimensional vector space. The proof uses the ellipticity of $\mathcal{M}$ to split the space of sections into a finite dimensional piece corresponding to the kernel of an elliptic operator (where Lojasiewicz's inequality holds) and an infinite dimensional piece orthogonal to the kernel (where Simon uses elliptic estimates).

In \cite{radethesis} and \cite{Rade92}, Rade extends this estimate to the case of the Yang-Mills functional on $2$ and $3$ dimensional manifolds. Simon's result does not hold \emph{a priori} since the gradient of the Yang-Mills functional is not an elliptic operator, however Rade uses a Coulomb gauge theorem to show that after the action of the gauge group one can restrict to a subspace where the Hessian is an elliptic operator, and then prove the result directly, following Simon's technique. 

In this section we prove Theorem \ref{thm:generalloja}, which is a Higgs bundle version of Simon's estimate for the functional $\QH$ defined below, and it is then shown that Proposition \ref{lem:loja} follows from Theorem \ref{thm:generalloja}. Many aspects of the proof of Theorem \ref{thm:generalloja} are more general than just the case of Higgs bundles considered in this paper, and can be extended to functionals on other spaces, such as the case of quiver bundles over Riemann surfaces (for which an analog of Hitchin and Simpson's theorem was proven in \cite{AlvarezGarcia03}). With this in mind, when possible the results are given in more general terms. 

For notation, let $M$ denote the affine Hilbert space $\left( T^* \mathcal{A} \right)_{H^1}$ and let $\mathcal{G}_{H^2}$ denote the completion of $\mathcal{G}$ in the $H^2$ norm. Note that Sobolev multiplication implies that $\mathcal{G}_{H^2}$ acts on $M$.

\begin{thm}\label{thm:generalloja}
Let $\rho: M \times \mathfrak{g} \rightarrow TM$ denote the infinitesimal action of $\mathcal{G}_{H^2}$ on $M$, and consider the functional $\QH:M \rightarrow \R$ defined by
\begin{equation*}
\QH(x) = \| \mu_1(x) \|^2 + \|\mu_2(x) \|^2 + \| \mu_3(x) \|^2
\end{equation*}
where $x$ denotes the point $(A'', \phi) \in \left( T^* \mathcal{A} \right)_{H^1}$.

Fix a critical point $x$ of $\QH$. Then there exists some $\varepsilon > 0$ (depending on $x$) and $\theta \in \left(0, \frac{1}{2} \right)$ such that the following inequality holds:
\begin{equation}\label{eqn:lojainequality}
\| \grad \QH (y) \|_{H^{-1}} \geq C \left| \QH(y) - \QH(x) \right|^{1-\theta} 
\end{equation}
whenever $\hone{x-y} < \varepsilon$.
\end{thm}

Assuming the result of the theorem, the proof of Proposition \ref{lem:loja} is as follows.

\begin{proof}[Proof of Proposition \ref{lem:loja}]

Choose a critical point $(A''_\infty, \phi_\infty) \in \mathcal{B}$ of $\YMH$, which is also a critical point of $\QH$. Note that $\left. \QH \right|_{\mathcal{B}} = \YMH$ and apply Theorem \ref{thm:generalloja} to show that there exists $\varepsilon > 0$ and $\theta \in \left(0, \frac{1}{2} \right)$  such that the following inequality holds for $(A, \phi) \in \mathcal{B}$ such that $\hone{(A'', \phi) - (A''_\infty, \phi_\infty)} < \varepsilon$.
\begin{equation*}
\| D_{(A'', \phi)}^* F_{(A'', \phi)} \|_{H^{-1}} \geq \left| \YMH(A, \phi) - \YMH(A_\infty, \phi_\infty) \right|^{1-\theta}
\end{equation*}
The inequality $\| D_{(A'', \phi)}^* F_{(A'', \phi)} \|_{H^{-1}} \leq \| D_{(A'', \phi)}^* F_{(A'', \phi)} \|_{L^2}$ completes the proof.
\end{proof}

The first step in the proof of Theorem \ref{thm:generalloja} is the following local description around a critical point.
\begin{prop}[Coulomb Gauge]\label{prop:coulombgauge}
Let $M$ be an affine Hilbert manifold with the action of a Hilbert Lie group $G$, and let $f:M \rightarrow \R$ be a $G$-invariant functional. Let $x\in M$ be a critical point of $f$ and  denote the Hessian of $f$ at the point $x\in M$ by $H_f(x): T_x M \rightarrow T_x M$. Let $\rho_x : \Lie(G) \rightarrow T_x M$ be the infinitesimal action of the group $G$ at the point $x \in M$ and suppose that the following operator is elliptic
\begin{equation}\label{eqn:Hessiangroupoperator}
H_f(x) + \rho_x \rho_x^* : T_x M \rightarrow T_x M
\end{equation}
Then there exists $\varepsilon > 0$ such that if $\|y-x\|<\varepsilon$ then there exists $u \in (\ker \rho_x)^\perp$ such that for $g = e^{-u}$ 
\begin{equation}
\rho_x^* \left( g \cdot y - x \right) = 0
\end{equation}

\end{prop}

This more general situation described above is related to the space of Higgs bundles in the following way. The functional $\QH$ is $\mathcal{G}_{H^2}$-invariant, and the Hessian is given (in the notation of Section \ref{sec:symplectic}) by the following formula
\begin{multline*}
\frac{1}{2} H_{\QH}(x) =  -I\rho_x \rho_x^* I X + I \delta \rho_x(*\mu_1(x))(X) -J\rho_x \rho_x^* J X \\
  + J \delta \rho_x(*\mu_2(x))(X) -K\rho_x \rho_x^* K X + K \delta \rho_x(*\mu_3(x))(X) 
\end{multline*}
From this description of the Hessian together with the description of the operator $\rho_x$ from Section \ref{sec:symplectic}  and complex structures $I,J,K$ from \cite{Hitchin87}, we see that $H_{\QH}(x) + \rho_x \rho_x^*$ is an elliptic operator on the tangent space 
$$T_x (T^* \mathcal{A})_{H^1} \cong H^1\left( \Omega^{0,1}(\End(E)) \oplus \Omega^{1,0}(\End(E)) \right)$$
and therefore the critical points of the functional $\QH$ on the space $T^* \mathcal{A}$ satisfy the conditions of Proposition \ref{prop:coulombgauge}. Since a critical point $(A'', \phi) \in \mathcal{B}$ of $\YMH$ is also a critical point of $\QH$, then the theorem applies at all critical points of $\YMH$.

The first step in the proof of Proposition \ref{prop:coulombgauge} is the following lemma.

\begin{lem}\label{lem:tangentdecomp}

At a critical point $x \in M$ of the functional $f$, $\im \rho_x$ is a closed subspace of $T_x M$ and the following decomposition holds
\begin{equation*}
T_x M \cong \ker \rho_x^* \oplus \im \rho_x
\end{equation*}

\end{lem}

The proof of this lemma in turn depends on the following lemmas.

\begin{lem}\label{lem:kernelL}
Let $L = H_f(x) + \rho_x \rho_x^*$. Then $\ker L = \ker H_f(x) \cap \ker \rho_x^*$.
\end{lem}

\begin{proof}[Proof of Lemma \ref{lem:kernelL}]
$f$ is $G$-invariant implies that $\im \rho_x \subseteq \ker H_f(x)$, and since the Hessian $H_f(x)$ is self-adjoint then $\im H_f(x) \subseteq \ker \rho_x^*$. Therefore $\im \rho_x \subseteq \left(\im H_f(x) \right)^\perp$, and so $\ker L \subseteq \ker H_f(x) \cap \ker \rho_x^*$. The inclusion $\ker H_f(x) \cap \ker \rho_x^* \subseteq \ker L$ follows from the definition of $L$.
\end{proof}

Using this lemma together with the fact that $L$ is elliptic and self-adjoint, we have the splitting
\begin{equation}\label{eqn:tangentsplitting}
T_x M = \ker L \oplus \im L^* \cong \ker L \oplus \im L \cong \left( \ker H_f(x) \cap \ker \rho_x^* \right) \oplus \im L
\end{equation} 

Next we need the following technical lemma.

\begin{lemma}\label{lem:closedsubspace}
Let $H$ be a closed Hilbert space with two linear subspaces $A,B \subset H$ that satisfy $A \subseteq B^\perp$ and $B \subseteq A^\perp$. If $H = A+B$, then $A$ and $B$ are closed subspaces, and $H = A \oplus B$.
\end{lemma}

\begin{proof}

The result follows from showing that $A = B^\perp$. Arguing by contradiction, suppose that $x \in B^\perp \setminus A$. Then $x = a + b$ for $a \in A$ and $b \in B$. Since $a \perp b$, then $\| x \|^2 = \|a\|^2 + \|b\|^2$. We can also write $a = x - b$, so $\|a\|^2 = \| x-b\|^2 = \|x\|^2 + \|b\|^2$, since $x \in B^\perp$ by assumption. Therefore $\| x\|^2 = \|x\|^2 + 2\|b\|^2$ and so $b=0$. This implies that $x \in A$, which is a contradiction.

Therefore $A = B^\perp$, and in particular $A$ is a closed subspace of $H$. Repeating the same argument shows that $B = A^\perp$. Since $H = A + B$ and $A,B$ are closed, orthogonal subspaces of $H$ then $H = A \oplus B$.
\end{proof}

\begin{lemma}\label{lem:decomposeimageL}

$\im L$ decomposes into a direct sum of closed subspaces
\begin{equation}\label{eqn:decomposeimageL}
\im L = \im H_f(x) \oplus \im \rho_x
\end{equation}

\end{lemma}

\begin{proof}[Proof of Lemma \ref{lem:decomposeimageL}]

Firstly note that 
\begin{equation}\label{eqn:imageinclusion}
\im L = \im (H_f(x) + \rho_x \rho_x^*) \subseteq \im H_f(x) + \im \rho_x \subseteq \left( \ker L \right)^\perp
\end{equation}
and since $L$ is elliptic, $\im L = \left( \ker L \right)^\perp$, and so all of the set inclusions in \eqref{eqn:imageinclusion} are equalities. Therefore $\im L$ is a closed Hilbert space such that $\im L = \im H_f(x) + \im \rho_x$. Recalling that $\im H_f(x) \perp \im \rho_x$ and applying Lemma \ref{lem:closedsubspace} completes the proof.
\end{proof}

\begin{proof}[Proof of Lemma \ref{lem:tangentdecomp}]

Applying Lemma \ref{lem:decomposeimageL} to the decomposition \eqref{eqn:tangentsplitting} shows that
\begin{equation}\label{eqn:furthertangentsplitting}
T_x M = \im \rho_x \oplus \im H_f(x) \oplus \left( \ker H_f(x) \cap \ker \rho_x^* \right)
\end{equation}
Since $\im H_f(x) \oplus \left( \ker H_f(x) \cap \ker \rho_x^* \right) \subseteq \ker \rho_x^*$, and $\im \rho_x \perp \ker \rho_x^*$ then applying Lemma \ref{lem:closedsubspace} to \eqref{eqn:furthertangentsplitting} gives us the decomposition
\begin{equation*}
T_x M = \im \rho_x \oplus \ker \rho_x^* \qedhere
\end{equation*}

\end{proof}

To complete the proof of Proposition \ref{prop:coulombgauge} we need the following description of a neighbourhood of the critical point $x$.
\begin{lem}\label{lem:localdiffeo}
The map $F : \left( \ker \rho_x \right)^\perp \times \ker \rho_x^* \rightarrow M$ given by 
\begin{equation}
F(u, X) = e^u \cdot (x + X)
\end{equation}
is a local diffeomorphism about the point $F(0, 0) = x$.
\end{lem}

\begin{proof}
$dF_{(0,0)} (\delta u, \delta X) = \rho_x(\delta u) + \delta X$. Since $\delta u \in \left( \ker \rho_x \right)^\perp$ and $\delta X \in \ker \rho_x^*$ then $dF_{(0,0)}$ is injective. By Lemma \ref{lem:tangentdecomp}, $T_x M \cong \ker \rho_x^* \oplus \im \rho_x$ and so $dF_{(0,0)}$ is an isomorphism. Applying the inverse function theorem completes the proof.
\end{proof}

\begin{proof}[Proof of Proposition \ref{prop:coulombgauge}]
Lemma \ref{lem:localdiffeo} shows that there exists $\varepsilon > 0$ such that if $\| y - x \| < \varepsilon$ then there exists $(u, X) \in \left( \ker \rho_x \right)^\perp \times \left( \ker \rho_x^* \right)$ such that $e^u \cdot (x + X) = y$. Re-arranging this gives us
\begin{equation}
X = e^{-u} \cdot y - x
\end{equation}
and since $X \in \ker \rho_x^*$, then setting $g = e^{-u}$ completes the proof.
\end{proof}

The function $\QH$ defined on $M$ satisfies the conditions of Proposition \ref{prop:coulombgauge}, and so at a critical point $x\in M$ we have the splitting $T_x M \cong \im \rho_x \oplus \ker \rho_x^*$. Using this decomposition of the tangent space, define projection operators $\Pi_{ker}$ and $\Pi_{im}$ denoting projection onto $\ker \rho_x^*$ and $\im \rho_x$ respectively. Since the inequality \eqref{eqn:lojainequality} is $\mathcal{G}_{H^2}$-invariant, then we can use Proposition \ref{prop:coulombgauge} to restrict attention to those points $y$ in a $\delta$-neighbourhood of $x$ such that $y-x \in \ker \rho_x^*$. Consider the functional $E : \ker \rho_x^* \rightarrow \R$ given by $E(b) = \QH(x+b) - \QH(b)$. The gradient of $E$ at the point $X \in \ker \rho_x^*$ is then given by $\grad E(b) = N(b) = \Pi_{ker} \grad \QH(x+b)$. Since the functional $\QH$ is analytic, then so is $E$ and hence $N$. The image of the Hessian of $\QH$ satisfies $\im H_{\QH}(x) \subseteq \ker \rho_x^*$, and so the derivative of $N$ at $b=0$ for $b' \in \ker \rho_x^*$ has the following expression
\begin{equation*}
dN_{b=0}(b') = \Pi_{ker} H_{\QH}(x)(b') = H_{\QH}(x)(b')
\end{equation*}
$H_{\QH}(x)$ is an elliptic operator $\ker \rho_x^* \rightarrow \ker \rho_x^*$ and so we can decompose $\ker \rho_x^*$ into closed subspaces
\begin{equation*}
\ker \rho_x^* \cong \left( \ker H_{\QH}(x) \cap \ker \rho_x^* \right) \oplus \im H_{\QH}(x)
\end{equation*}
For notation, write $K_0 = \left( \ker H_{\QH}(x) \cap \ker \rho_x^* \right)$ and decompose $\ker \rho_x^* \cap H^s \cong K_0 \oplus K_\pm^s$. Denote the norm on $K_0$ by $\| \cdot \|_{K_0}$ and note that since $K_0$ is the kernel of an elliptic operator then it is finite dimensional and all norms on $K_0$ are equivalent. For any $b \in \ker \rho_x^* \cap H^1$ write $b = b_0 + b_\pm$ with $b_0 \in K_0$ and $b_\pm \in K_\pm^1$. Since $H_{\QH}(x)$ is the derivative of $N: \ker \rho_x^* \rightarrow \ker \rho_x^*$ and $H_{\QH}(x)$ is an injective operator $K_\pm^1 \rightarrow K_\pm^{-1}$ (and so an isomorphism onto its image), then an application of the implicit function theorem gives us the following lemma.

\begin{lemma}\label{lem:implicitfunction}

There exists $\varepsilon > 0$ and $\delta > 0$, and a map $\ell: B_{\varepsilon} K_0 \rightarrow B_{\delta} K_\pm^1$, such that for any $b_0 \in K_0$ satisfying $\| b_0 \|_{K_0} < \varepsilon$, we have that $N(b) \in K_0$ if and only if $b = b_0 + \ell (b_0)$. Moreover, since the function $N$ is analytic then so is $\ell$.

\end{lemma}
Given any $b \in \ker \rho_x^*$ we can use Lemma \ref{lem:implicitfunction} to decompose $b = b_0 + \ell(b_0) + b_\pm$, where $b_0$ is the projection of $b$ onto the subspace $K_0 \subseteq T_x M$ and $b_\pm = b - b_0 - \ell(b_0)$

\begin{lem}
There exists $\varepsilon > 0$ and $\delta > 0$ such that for $\| b \|_{H^1} < \varepsilon$ the following inequalities hold
\begin{equation*}
\| b_0 \|_{K_0} \leq c \| b \|_{H^1}, \quad
\| \ell(b_0) \|_{H^1} \leq c \| b \|_{H^1}, \quad
\| b_\pm \|_{H^1} \leq c \| b \|_{H^1} 
\end{equation*}
\end{lem}

\begin{proof} 
$K_0 \subseteq T_x M$ is finite-dimensional, therefore all norms on $K_0$ are equivalent and there exists $c$ such that $\| b_0 \|_{K_0} \leq c \| b_0 \|_{H^1}$. Also, $b_0 \perp (\ell(b_0) + b_\pm )$ implies that $\| b_0 \|_{H^1} + \| \ell(b_0) + b_\pm \|_{H^1} = \| b \|_{H^1}$, so therefore $\| b_0 \|_{H^1} \leq \| b \|_{H^1}$. Since $\ell : B_{\varepsilon} K_0 \rightarrow B_{\delta} K_\pm$ is smooth and has a finite dimensional domain then for some $k$ we have $\| \ell(b_0) \|_{H^1} \leq k \| b_0 \|_{K_0} \leq (c-1) \| b \|_{H^1}$. Therefore there exists a constant $c$ such that $\| b_\pm \|_{H^1} = \| b - b_0 - \ell(b_0) \|_{H^1} \leq \| b \|_{H^1} + \| \ell(b_0) \|_{H^1} + \| b_0 \|_{H^1} \leq c \| b \|_{H^1}$.
\end{proof}

Denote the completion of $T_x M$ in the $H^s$ norm by $(T_x M)_{H^s}$. Define $g: K_0 \rightarrow \R$ by $g(b_0) = E(b_0 + \ell(b_0) )$ and note that since $E$ and $\ell$ are real analytic then $g$ is real analytic. Now we can split $N(b)$ into the following parts
\begin{align}\label{eqn:splitM}
N(b) & = N(b_0 + \ell(b_0) + b_\pm) \nonumber \\
 & = \nabla g(b_0) - N(b_0 + \ell(b_0) ) + N(b_0 + \ell(b_0) + b_\pm) \nonumber \\
 & = \nabla g(b_0) + \int_0^1 dN(b_0 + \ell(b_0) + sb_\pm)(b_\pm) \, ds \nonumber \\
 & = \nabla g(b_0) + H_{\QH}(x)(b_\pm) + L_1(b_\pm)
\end{align}
where $L_1 : (T_x M)_{H^1} \rightarrow (T_x M)_{H^{-1}}$ is defined by $L_1(a) = \int_0^1 dN(b_0 + \ell(b_0) + sb_\pm)(a)  - dN(0)(a) \, ds$.

\begin{cla}\label{cla:L1bounded}
\begin{equation*}
\| L_1(b_\pm) \|_{H^{-1}} \leq c \| b \|_{H^1} \| b_\pm \|_{H^1}
\end{equation*}
\end{cla}

\begin{proof}
For $b \in T_x M$ define $h_s(b) = dN(b_0 + \ell(b_0) + s b_\pm) - dN(0)$. Since $N$ is analytic then $h_s$ is also analytic, and together with the fact that $h_s(0) = 0$ then there exists $\varepsilon > 0$ and some constant $c(s)$ depending on $s \in [0,1]$ such that whenever $\| b \|_{H^1} < \varepsilon$ we have the following inequality 
\begin{align*}
\| h_s(b_0 + \ell(b_0) + b_\pm) \|_{H^{-1}} & \leq c(s) \| b_0 + \ell(b_0) + b_\pm \|_{H^1} \\
 \Rightarrow \| h_s(b_0 + \ell(b_0) + b_\pm) \|_{H^{-1}} & \leq C \| b \|_{H^1} \\
\Rightarrow \, \| dN(b_0 + \ell(b_0) + s b_\pm) - dN(0) \|_{H^{-1}} & \leq C \| b \|_{H^1}  
\end{align*}
Therefore $\| L_1(b_\pm) \|_{H^{-1}} \leq C \| b \|_{H^1} \| b_\pm \|_{H^1}$ whenever $\| b \|_{H^1} < \varepsilon$.
\end{proof}

\begin{lemma}
The following inequality holds whenever $\| b \|_{H^1} < \varepsilon$
\begin{equation}\label{eqn:Nboundedbelow}
\| N(b) \|_{H^{-1}} \geq c \left( \| \nabla g(b_0) \|_{K_0} + \| b_\pm \|_{H^1} \right) 
\end{equation}
\end{lemma}

\begin{proof}

By Lemma \ref{lem:implicitfunction}, $\nabla g(b_0) \in K_0$ and $H_{\QH}(b_\pm) \in K_\pm^{-1}$. Therefore $\nabla g(b_0) \perp H_{\QH}(b_\pm)$, which together with \eqref{eqn:splitM} implies that
\begin{equation}\label{eqn:splitinequalityN}
\| N(b) \|_{H^{-1}} \geq \| \nabla g(b_0) \|_{H^{-1}} + \| H_{\QH}(x)(b_\pm) \|_{H^{-1}} - \| L_1(b_\pm) \|_{H^{-1}}
\end{equation}
Since all norms on $K_0$ are equivalent then $\| \nabla g(b_0) \|_{H^{-1}} \geq c \| \nabla g(b_0) \|_{K_0}$ for some constant $c$. $K_\pm^1$ is orthogonal to $\ker H_{\QH}$ and $H_{\QH}$ is elliptic, therefore $\| H_{\QH}(x)(b_\pm) \|_{H^{-1}} \geq c \| b_\pm \|_{H^1}$ for some constant $c$. Together with \eqref{eqn:splitinequalityN} and Claim \ref{cla:L1bounded} this completes the proof.
\end{proof}

We can decompose the functional $E$ in the following way
\begin{align}\label{eqn:decomposeE}
E(b) & = g(b_0) + E(b_0 + \ell(b_0) + b_\pm) - E(b_0 + \ell(b_0)) \nonumber \\
 & = g(b_0) + \int_0^1 \left< N(b_0 + \ell(b_0) + s b_\pm), b_\pm \right>  ds \nonumber \\
 & = g(b_0) + \left< N(b_0 + \ell(b_0) ), b_\pm \right> \nonumber \\
 & \, \, \, + \int_0^1 \left< N(b_0 + \ell(b_0) + s b_\pm) - N(b_0 + \ell(b_0)), b_\pm \right> ds \nonumber  \\
 & = g(b_0) + \left< \nabla g(b_0), b_\pm \right> \nonumber \\
 & \, \, \, + \int_0^1 \int_0^1 \left< dN(b_0 + \ell(b_0) + s t b_\pm) (sb_\pm), b_\pm \right> ds \, dt \nonumber \\
 & = g(b_0) + \left< \nabla g(b_0), b_\pm \right> + \left< H_{\QH}(x)(b_\pm), b_\pm \right> + \left< L_2(b_\pm), b_\pm \right> 
\end{align} 
where $L_2 : (T_x M)_{H^1} \rightarrow (T_x M)_{H^{-1}}$ is defined  by
\begin{equation*}
L_2(b_\pm) = \int_0^1 \int_0^1 \left( dN(b_0 + \ell(b_0) + s t b_\pm)(s b_\pm) - H_{\QH}(b_\pm) \right)  ds \, dt
\end{equation*}

\begin{lem}
The following holds whenever $\| b \|_{H^1} < \varepsilon$
\begin{equation}\label{eqn:boundonE}
\left| E(b) \right| \leq \left| g(b_0) \right| + C \| b_\pm \|_{H^1}^2
\end{equation}
\end{lem} 

\begin{proof}

Following the same proof as Claim \ref{cla:L1bounded}, we have that whenever $\hone{b} < \varepsilon$
\begin{equation*}
\| L_2(b_\pm) \|_{H^{-1}} \leq c \| b \|_{H^1} \| b_\pm \|_{H^1}
\end{equation*}
Since $\nabla g(b_0) \in K_0$ and $b_\pm \in K_\pm$ then $\left< \nabla g(b_0), b_\pm \right> = 0$. $H_{\QH}(x)$ is elliptic and injective on $K_\pm^1$, therefore $\left< H_{\QH}(x)(b_\pm), b_\pm \right> \leq c \| b_\pm \|_{H^1}^2$. Applying these results to \eqref{eqn:decomposeE} completes the proof.
\end{proof}
 
 \begin{proof}[Proof of Theorem \ref{thm:generalloja}]
 
Since the inequality \eqref{eqn:lojainequality} is $\mathcal{G}_{H^2}$-invariant, then we can use Proposition \ref{prop:coulombgauge} to restrict to those points $y$ in a $\delta$-neighbourhood of $x$ such that $y-x \in \ker \rho_x^*$.  Therefore it is sufficient to prove that $\left| E(b) \right|^{1-\theta} \leq K \| N(b) \|_{H^{-1}}$ for some constant $K$ and $b \in \ker \rho_x^* \cap H^1$ with $\| b \|_{H^1} < \varepsilon$. Since $g:K_0 \rightarrow \R$ is a real analytic function on a finite-dimensional space, from results of Lojasiewicz in \cite{Lojasiewicz84} there exists $\theta \in (0, \frac{1}{2})$ such that
\begin{equation}\label{eqn:originallojasiewiczinequality}
\left| g(b_0) \right|^{1-\theta} \leq c \| \nabla g(b_0) \|_{K_0}
\end{equation}
Applying this to equation \eqref{eqn:boundonE} and using \eqref{eqn:Nboundedbelow} gives us
\begin{align*}
\left| E(b) \right|^{1-\theta} & \leq \left|  g(b_0) \right|^{1-\theta} + c \| b_\pm \|_{H^1}^{2(1-\theta)} \\
 & \leq c \left( \| \nabla g(b_0) \|_{K_0} + \| b_\pm \|_{H^1} \right) \\
 & \leq K \| N(b) \|_{H^{-1}}
\end{align*}
for any $b$ in $\ker \rho_x^* \cap H^1$ with $\| b \|_{H^1} < \varepsilon$.  
\end{proof}

\subsection{An interior estimate}\label{subsec:energy}

The purpose of this section is to prove Proposition \ref{lem:energy}, which provides an estimate relating the $H^k$ and $L^2$ norms of a tangent vector to the gradient flow of $\YMH$. The relationship between the $H^1$ and $L^2$ norms of a tangent vector to the Yang-Mills flow was proved in \cite{Rade92}, here we extend these results to derive estimates on higher derivatives of the gradient of the functional $\YMH$ on the space $\mathcal{B}$. 

Recall that the proof of Proposition \ref{lem:convmodgauge} shows that $\| F_A + [\phi, \phi^*] \|_{H^k}$ is bounded for all $k$. For fixed $(A_0'', \phi_0)$ define the bounded complex gauge orbit 
\begin{equation}
\mathcal{O}_C^k = \left\{ (A'', \phi) \in \mathcal{G}^\C \cdot(A_0'', \phi_0) : \| F_A + [\phi, \phi^*] \|_{H^k} < C \right\}
\end{equation}
As noted in Subsection \ref{subsec:compactness}, the proof of Lemma \ref{lem:L4compactness} shows that there exists a constant $K$ such that $\| F_A \|_{H^k} < K$ on $\mathcal{O}_C^k$. For $(A'', \phi) \in \mathcal{O}_C^k$ consider the initial value problem
\begin{equation}\label{eqn:initialvalueproblem}
\left\{ \begin{matrix} \frac{\partial \psi}{\partial t} + \nabla_A^* \nabla_A \psi = \sigma \\ \psi(0, \cdot) = 0 \end{matrix} \right.
\end{equation}
for some $\sigma \in H^s$ where $s \in [-k-2, k+2]$. Equations (11.3), (11.4) and Proposition A in the Appendix of \cite{Rade92} show that the following estimates hold on the time interval $[0, 2t_0]$
\begin{align}
\| \psi \|_{L^2([0, 2t_0], H^k)} & \leq c_K \| \sigma \|_{L^1([0, 2t_0], H^{k-1})} \label{eqn:initialvalueestimateone} \\
\| \psi \|_{L^2([0, 2t_0], H^k)} & \leq c_K t_0^{\frac{1}{4}} \| \sigma \|_{L^2([0, 2t_0], H^{k-\frac{3}{2}})} \label{eqn:initialvalueestimatetwo}
\end{align}
Moreover the constant $c_K$ only depends on $K$. For notation, in the following we use $H^{0,k}$ to denote the space $L^2 \left([0, 2t_0], H^k \right)$. Now let $G(t) = \left( \frac{\partial A}{\partial t}, \frac{\partial \phi}{\partial t} \right)$ and using the notation of Section \ref{sec:symplectic}, for any moment map on a symplectic manifold we have the following formula for $\frac{\partial G}{\partial t}$ along the downwards gradient flow of $\| \mu \|^2$
\begin{align*}
G(t) & = \frac{\partial x}{\partial t} = -I \rho_x(*\mu) \\
\Rightarrow \, \, \frac{\partial G}{\partial t} & = -I \delta \rho_x(*\mu)\left( \frac{\partial x}{\partial t} \right) - I \rho_x \left( *d \mu\left( \frac{\partial x}{\partial t} \right) \right) \\
 & = - I \delta \rho_x(*\mu)(G(t)) + I \rho_x \rho_x^* I (G(t))
\end{align*}
where in the last step we identify $*d\mu = -\rho_x^* I$ as described in Section \ref{sec:symplectic}. Note also that equation \eqref{eqn:infadjoint} shows that $\rho_x^* G = 0$. Therefore we have the equation
\begin{equation}\label{eqn;timederivgradient}
\frac{\partial G}{\partial t} + \rho_x \rho_x^* (G) - I \rho_x \rho_x^* I (G) = -I \delta \rho_x(*\mu)(G)
\end{equation}
which for the case of $\mu = F_A + [\phi, \phi^*]$ reduces to
\begin{equation}\label{eqn:timederivhiggsgradient}
\frac{\partial G}{\partial t} + \Delta_{(A'', \phi)} G = \left\{ F_A + [\phi, \phi^*], G \right\}
\end{equation}
where the operator $\left\{ \cdot , \cdot \right\}$ denotes various different intrinsically defined multilinear operators. The Weitzenb\"ock formula of Simpson (Lemma 7.2.1 of \cite{simpsonthesis}) states that for a $k$-form $\alpha$ with values in $E$
\begin{multline}\label{eqn:simpsonweitz}
\nabla_A^* \nabla_A \alpha = \Delta_{(A'', \phi)} \alpha + \left\{ F_A, \alpha \right\} + \left\{ (\phi + \phi^*), (\phi + \phi^*), \alpha \right\} \\ + \left\{ \nabla_A \phi, \alpha \right\} + \left\{ R, \alpha \right\},
\end{multline}
where $R$ refers to the Riemannian curvature of $X$. Substituting this formula into \eqref{eqn:timederivhiggsgradient} gives the following expression
\begin{equation}\label{eqn:initialvaluegradient}
\frac{\partial G}{\partial t} + \nabla_A^* \nabla_A G = \left\{ F_A + [\phi, \phi^*], G \right\} + \left\{ \nabla_A \phi, G \right\} + \left\{ R, G \right\}
\end{equation}
Now let $(A'', \phi) = (A_\infty'', \phi_\infty) + a$ (where $a \in \Omega^1(\End(E))$, and note that the assumption of Proposition \ref{lem:energy} is that $\| a \|_{H^k} < \varepsilon_1$. Then we have the following equation
\begin{multline}
\frac{\partial G}{\partial t} + \nabla_\infty^* \nabla_\infty G = \left\{ F_A + [\phi, \phi^*], G \right\} + \left\{ \nabla_A \phi, G \right\} + \left\{ R, G \right\} \\
 + \left\{ a, \nabla_\infty G \right\} + \left\{ \nabla_\infty a, G \right\} + \left\{ a, a, G \right\}
\end{multline}
Multiplying both sides by a smooth cut-off function $\eta(t)$ with $\eta = 0$ on $[0, \frac{1}{2}t_0]$ and $\eta = 1$ on $[t_0, 2t_0]$ gives the equation
\begin{multline}
\frac{\partial (\eta G)}{\partial t} + \nabla_\infty^* \nabla_\infty (\eta G) = \left\{ F_A + [\phi, \phi^*], \eta G \right\} + \left\{ \nabla_A \phi, \eta G \right\} + \left\{ R, \eta G \right\} \\
 + \left\{ a, \nabla_\infty ( \eta G ) \right\} + \left\{ \nabla_\infty a, \eta G \right\} + \left\{ a, a, \eta G \right\} + \frac{\partial \eta}{\partial t} G 
\end{multline}
The existence of a solution to the gradient flow equations \eqref{eqn:gradfloweqns} shows that $\eta G$ is a solution to the initial value problem \eqref{eqn:initialvalueproblem}. Therefore, following the method of \cite{Rade92} p156 (see also \cite{Woodward06} p30 for more details), the estimates \eqref{eqn:initialvalueestimateone} and \eqref{eqn:initialvalueestimatetwo} show that
\begin{multline*}
\| \eta G \|_{H^{0,k}} \leq C t_0^{\frac{1}{4}} \| \left\{ F_A + [\phi, \phi^*], \eta G \right\} + \left\{ \nabla_A \phi, \eta G \right\} + \left\{ R, \eta G \right\} \|_{H^{0,k-\frac{3}{2}}} \\
 \quad \quad \quad + C t_0^{\frac{1}{4}} \| \left\{ a, \nabla_\infty ( \eta G ) \right\} + \left\{ \nabla_\infty a, \eta G \right\} + \left\{ a, a, \eta G \right\} \|_{H^{0,k-\frac{3}{2}}} \\
 + C \| \frac{\partial \eta}{\partial t} G \|_{L^1([0, 2 t_0], H^{k-1})}
\end{multline*}
Therefore Sobolev multiplication theorems as used in \cite{Rade92} show that for $G \in H^k$, $F_A + [\phi, \phi^*] \in H^{k-1}$, $\nabla_A \phi \in H^{k-1}$, $R$ is smooth and $a \in H^k$, the following two estimates hold
\begin{align*}
& \| \left\{ F_A + [\phi, \phi^*], \eta G \right\} + \left\{ \nabla_A \phi, \eta G \right\} + \left\{ R, \eta G \right\} \|_{H^{0,k-\frac{3}{2}}} \leq C \| \eta G \|_{H^{0,k}} \\
\quad & \| \left\{ a, \nabla_\infty ( \eta G ) \right\} + \left\{ \nabla_\infty a, \eta G \right\} + \left\{ a, a, \eta G \right\} \|_{H^{0, k-\frac{3}{2}}} \\
  & \quad \quad \quad \leq C \| \eta G \|_{H^{0,k}} + \varepsilon_1 \| \nabla_\infty (\eta G) \|_{H^{0,k-1}} \leq C \| \eta G \|_{H^{0,k}}
\end{align*}

Therefore $\| \eta G \|_{H^{0, k}} \leq C t_0^{\frac{1}{4}} \| \eta G \|_{H^{0, k}} + C \| \frac{\partial \eta}{\partial t} G \|_{L^1([0, 2t_0], H^{k-1})}$, and so when $t_0$ is small we have $\| \eta G \|_{H^{0, k}} \leq C' \| \frac{\partial \eta }{\partial t} G \|_{L^1([0, 2t_0], H^{k-1})}$. Therefore
\begin{align}
\begin{split}
\| G \|_{L^1([t_0, 2t_0], H^k)} & \leq C t_0^{\frac{1}{2}} \| \eta G \|_{L^2([0, 2t_0], H^k)} \\
&  \leq Ct_0^{\frac{1}{2}} \| \frac{\partial \eta}{\partial t} G \|_{L^1([0, 2t_0], H^{k-1})} \leq C t_0^{-\frac{1}{2}} \| G \|_{L^1([0, 2t_0], H^{k-1})}
\end{split}
\end{align}
Dividing the interval $[T,S]$ into subintervals of length $t_0$ and applying this estimate on each sub-interval shows that (for $t_0 \leq \frac{1}{k}$)
\begin{align*}
\int_{T+1}^S \| G \|_{H^k} \, dt \leq \int_{T+ k t_0}^S \| G \|_{H^k} \, dt & \leq  C \int_{T + (k-1)t_0}^S \| G \|_{H^{k-1}} \, dt \\
 & \leq C \int_T^S \| G \|_{L^2} \, dt
\end{align*}
which completes the proof of Proposition \ref{lem:energy}.

\section{Algebraic and analytic stratifications}\label{sec:stratification}

In order to set the notation we first recall the main points of the Harder-Narasimhan filtration for Higgs bundles from \cite{HauselThaddeus04}. Given a filtration $E_0 \subset E_1 \subset \cdots \subset E_r = E$ of $E$ by $\phi$-invariant holomorphic sub-bundles, let $F_i = E_i / E_{i-1}$ and let $\phi_i \in \Omega^{1,0}(\End(F_i)$ be the induced Higgs field. The filtration is called a \emph{$\phi$-invariant Harder-Narasimhan filtration} if the pairs $(F_i, \phi_i)$ are semi-stable, and the slope $\frac{\deg(F_i)}{\rank(F_i)}$ is strictly decreasing in $i$. For a rank $n$ bundle, the \emph{type} of the Harder-Narasimhan filtration is the $n$-tuple $\mu = (\mu_1, \ldots, \mu_n)$, where the first $\rank(F_1)$ terms are $\frac{\deg(F_1)}{\rank(F_1)}$, the next $\rank(F_2)$ terms are $\frac{\deg(F_2)}{\rank(F_2)}$ and so on. Let $\mathcal{B}_\mu$ denote the space of Higgs pairs which have a $\phi$-invariant Harder-Narasimhan filtration of type $\mu$. As shown in \cite[Section 7]{HauselThaddeus04}, each Higgs pair possesses a unique Harder-Narasimhan filtration, the space $\displaystyle{\mathcal{B} = \bigcup_\mu \mathcal{B}_\mu}$ is stratified by these subsets, and the strata satisfy the closure condition $\displaystyle{\overline{\mathcal{B}_\mu} \subseteq \bigcup_{\nu \geq \mu} \mathcal{B}_\nu}$, where we use the usual partial ordering on Harder-Narasimhan types (cf Section 7 of \cite{AtiyahBott83} for holomorphic bundles, or Section 7 of \cite{HauselThaddeus04} for Higgs bundles). 

At a critical point $(A'', \phi)$ of $\YMH$, the bundle $E$ splits into $\phi$-invariant holomorphic sub-bundles and the goal of this section is to show that the \emph{algebraic stratification} by the type of the Harder-Narasimhan filtration is equivalent to the \emph{analytic stratification} by the gradient flow described in the previous section, where the equivalence is by the type of the splitting of $E$ into $\phi$-invariant holomorphic sub-bundles at the critical points of the functional $\YMH$.

In order to describe the analytic stratification of $\mathcal{B}$ using the results of Section \ref{sec:analysis}, firstly recall the critical point equations
\begin{align}
d_A''*(\muone) & = 0 \label{eqn:criteqn1} \\
\left[ \phi, *(\muone) \right] & = 0 \label{eqn:criteqn2}
\end{align}
Equation (\ref{eqn:criteqn1}) shows that for a non-minimal critical point $(A'', \phi)$ the bundle $E$ splits holomorphically into sub-bundles (see for example Theorem 3.1 in \cite{FreedUhlenbeck84} for the Yang-Mills functional), and equation \eqref{eqn:criteqn2} shows that the holomorphic sub-bundles are $\phi$-invariant. Therefore the space $\mathcal{B}_{crit}$ of non-minimal critical sets can be stratified by the Harder-Narasimhan type of each $\phi$-invariant holomorphic splitting $\displaystyle{\mathcal{B}_{crit} = \bigcup_\mu \eta_\mu}$. Given a Higgs pair $(A'', \phi)$ let $r(A'', \phi)$ denote the limit of the gradient flow with initial conditions $(A'', \phi)$ as defined in Section \ref{sec:analysis}. Define the \emph{analytic stratum} associated to each critical set by
\begin{equation}
\mathcal{C}_\mu = \left\{ (A, \phi) \in \mathcal{B} : r(A, \phi) \in \eta_\mu \right\}.
\end{equation}
Then Theorem \ref{thm:convergence} shows that $\mathcal{B}$ is stratified by the sets $\left\{ \mathcal{C}_\mu \right\}$ in the sense of Proposition 1.19 (1)-(4) of \cite{AtiyahBott83} (statement (5) of Proposition 1.19 in \cite{AtiyahBott83}, that the strata have well-defined codimension, cannot be true for $\left\{ \mathcal{C}_\mu \right\}$ because the dimension of the negative eigenspace of the Hessian of $\YMH$ is not constant). Moreover, each stratum $\mathcal{C}_\mu$ retracts $\mathcal{G}$-equivariantly onto the corresponding critical set $\eta_\mu$ with the retraction defined by the gradient flow. The main theorem to be proved in this section is the following.

\begin{thm}\label{prop:equivstratification}
The algebraic stratification by the $\phi$-invariant Harder-Narasimhan type $\left\{ \mathcal{B}_\mu \right\}$  is the same as the analytic stratification $\left\{ \mathcal{C}_\mu \right\}$ by the gradient flow of $\YMH$.
\end{thm}

The proof of the theorem relies on the following results. Let $\mathfrak{g}$ denote the Lie algebra of the structure group of $E$ (which will be $\mathfrak{u}(n)$ or $\mathfrak{su}(n)$ in our case) and note that the following analog of Proposition 8.22 from \cite{AtiyahBott83} also holds for the functional $\YMH$. 

For a pair $(A'', \phi)$ of type $\mu$ and a convex invariant function $h:\mathfrak{g} \rightarrow \R$, let $H(A'', \phi) = \inf \int_M h(*(F_A + [\phi, \phi^*]))$, where the infimum runs over all pairs $(A'', \phi) \in \mathcal{B}_\mu$. Also, if $\mu$ can be written as $\mu = (\lambda_1, \ldots, \lambda_n)$, let $\Lambda_\mu$ be the diagonal matrix with entries $-2\pi i \lambda_i$.
\begin{prop}\label{prop:convexdesc}
$(A'', \phi)$ is of type $\mu$ iff $H(A'', \phi) = h(\Lambda_\mu)$ for all convex invariant $h$. Moreover, $\left< \grad H, \grad \YMH \right> \geq 0$.
\end{prop}

The details are the same as those in Section 8 of \cite{AtiyahBott83} for the case of holomorphic bundles, and so the proof is omitted.

\begin{cla}\label{cla:lowerboundlambda}
If ${\mathcal{B}_\mu} \cap \mathcal{C}_\lambda$ is non-empty then $\lambda \geq \mu$.
\end{cla}

\begin{proof}[Proof of Claim \ref{cla:lowerboundlambda}]
Let $(A'', \phi) \in \mathcal{B}_\mu \cap \mathcal{C}_\lambda$ for $\lambda \neq \mu$. The proof of Proposition \ref{lem:existence} shows that finite-time gradient flow is equivalent to the action of an element of $\mathcal{G}^\C$. Therefore we can find $\{ g_j \} \subset \mathcal{G}^\C$ such that $g_j \cdot (A'', \phi) \rightarrow (A_\infty, \phi_\infty)$, and since $(A'', \phi) \in \mathcal{C}_\lambda$ then $(A_\infty, \phi_\infty)$ is of type $\lambda$. $(A'', \phi) \in {\mathcal{B}_\mu}$ and $\mathcal{G}^\C$ preserves ${\mathcal{B}_\mu}$, therefore by Proposition \ref{prop:convexdesc}, we have $\YMH \left(g_j \cdot(A'', \phi) \right) \geq \YMH(\Lambda_\mu)$ for all $j$. Therefore $\YMH(\Lambda_\lambda) = \YMH(A_\infty, \phi_\infty) \geq \YMH(\Lambda_\mu)$ also.
\end{proof} 

\begin{lem}\label{lem:localstrata}
For each Harder-Narasimhan type $\mu$, there exists a neighbourhood $V_\mu$ of $\eta_\mu$ such that $\mathcal{B}_\mu \cap V_\mu \subseteq \mathcal{C}_\mu$.
\end{lem}

\begin{proof}[Proof of Lemma \ref{lem:localstrata}]
The proof follows that of \cite{Daskal92} Proposition 4.12 for the Yang-Mills functional. Since the vector bundle $E$ has finite rank then the set $\{ \lambda_i \}$ such that $\YMH(\Lambda_{\lambda_i}) = \YMH(\Lambda_\mu)$ is finite. Choose $\varepsilon > 0$ such that the only critical sets $\eta_\lambda$ intersecting $U_\mu = \YMH^{-1} \left( \YMH(\Lambda_\mu) - \varepsilon, \YMH(\Lambda_\mu) + \varepsilon \right)$ are those for which $\lambda \in \{ \lambda_i \}$. By Claim \ref{cla:lowerboundlambda} we can restrict attention to those $\lambda$ for which $\lambda \geq \mu$. For each $\lambda \in \{ \lambda_i \}_{\lambda_i \geq \mu}$ choose a convex invariant functional $f_\lambda$ such that $f_\lambda(\lambda) > f_\lambda(\mu)$, and let $e_\lambda = \frac{1}{2}\left( f_\lambda(\lambda) - f_\lambda(\mu) \right)$. Define the sets
\begin{equation*}
V_\lambda = U_\mu \cap F_\lambda^{-1} \left( f_\lambda(\mu) - e_\lambda, f_\lambda(\mu) + e_\lambda \right)
\end{equation*}
and note that $\eta_\mu \subset V_\lambda$ for each $\lambda \in \{ \lambda_i \}_{\lambda_i \geq \mu}$. Suppose that $(A'', \phi) \in V_\lambda \cap {\mathcal{B}_\mu} \cap \mathcal{C}_\lambda$ and let $(A_\infty'', \phi_\infty)$ denote the limit of $(A'', \phi)$ under the gradient flow of $\YMH$. Therefore $f_\lambda(\lambda) = F_\lambda(A_\infty'', \phi_\infty) \leq F_\lambda(A'', \phi) < f_\lambda(\mu) + e_\lambda$, since $(A'', \phi) \in V_\lambda$. We then have
\begin{equation*}
f_\lambda(\lambda) < f_\lambda(\mu) + \frac{1}{2}\left( f_\lambda(\lambda) - f_\lambda(\mu) \right) = \frac{1}{2} \left( f_\lambda(\mu) + f_\lambda(\lambda) \right) < f_\lambda(\lambda)
\end{equation*}
a contradiction. Therefore $V_\lambda \cap {\mathcal{B}_\mu} \cap \mathcal{C}_\lambda = \emptyset$, and setting $\displaystyle{V_\mu = \bigcap_{\lambda \in \{ \lambda_i \}_{\lambda_i \geq \mu}} V_\lambda}$ completes the proof.
\end{proof}

\begin{lemma}\label{lem:existsmetric}
Let $(A'', \phi)$ be a Higgs pair of $\phi$-invariant Harder-Narasimhan type $\mu$, and let $\{ f_i \}$ be a finite collection of convex invariant functions as defined in Proposition \ref{prop:convexdesc}. Then for any $\varepsilon > 0$ there exists a metric $H$ on $E$ such that $f_i(F_H + [\phi_H,  \phi_H^*] - \Lambda_\mu) < \varepsilon$ for all $i$.
\end{lemma}

\begin{proof}
Theorem 1 of \cite{Simpson88} shows that the result holds if $(A'', \phi)$ is stable. For the case where $(A'', \phi)$ is semistable, the proof of Theorem 1 on p895 of \cite{Simpson88} shows that the functional $M(K, H_t)$ is bounded below, and $\frac{\partial}{\partial t} M(K, H_t) = - \| F_{H_t} + [\phi_{H_t}, \phi_{H_t}^*] - \mu \cdot \id \|_{L^2}$. Therefore $\int_t^{t+1} \| F_{H_t} + [\phi_{H_t}, \phi_{H_t}^*] - \mu \cdot \id \|_{L^2} \, dt \rightarrow 0$ as $t \rightarrow \infty$. Equation \eqref{eqn:normcurvatureheat} shows that Theorem \ref{thm:moserharnack} applies to the function $\left| F_{H_t} + [\phi_{H_t}, \phi_{H_t}^*] - \mu \cdot \id \right|$ and so $\| F_{H_t} + [\phi_{H_t}, \phi_{H_t}^*] - \mu \cdot \id \|_{C^0} \rightarrow 0$, which is enough to prove Lemma \ref{lem:existsmetric} for the semistable case (cf Corollary 25 of \cite{Donaldson85} for the case of the Yang-Mills functional on a K\"ahler surface).

For the case of a Higgs pair of general $\phi$-invariant Harder-Narasimhan type the result follows by induction on the length of the Harder-Narasimhan filtration, as in the proof of Theorem 3.10 in \cite{DaskalWentworth04} for the Yang-Mills functional.
\end{proof}

Applying this result to the functions $f_\lambda$ in the proof of Lemma \ref{lem:localstrata} gives
\begin{corollary}
Given $(A'', \phi) \in \mathcal{B}_\mu$ there exists $g_0 \in \mathcal{G}^\C$ such that $g_0 \cdot (A'', \phi) \in V_\mu$, and so $r(g_0 \cdot(A'', \phi)) \in \eta_\mu$.
\end{corollary}
The next lemma shows that if the $\mathcal{G}^\C$-orbit of $(A'', \phi)$ intersects $V_\mu$ then the gradient flow with initial conditions $(A'', \phi)$ converges to the critical set $\eta_\mu \subset V_\mu$.

\begin{lemma}\label{lem:convergeforallmetrics}
If there exists $g_0 \in \mathcal{G}^\C$ such that $r(g_0 \cdot (A'', \phi)) \in \eta_\mu$ then $r(g \cdot (A'', \phi)) \in \eta_\mu$ for all $g \in \mathcal{G}^\C$.
\end{lemma}

\begin{proof}
As noted in Section \ref{subsec:existence}, the action of an element $g \in \mathcal{G}^\C$ can be described up to $\mathcal{G}$-equivalence by changing the metric on $E$ by $H \mapsto H h$. Since the set $\eta_\mu$ is preserved by $\mathcal{G}$ and the gradient flow is $\mathcal{G}$-equivariant, then it is immediate that the lemma holds for all $g \in g_0 \cdot \mathcal{G}$, and so it is sufficient to show that the lemma is true for any Hermitian metric $H$ on $E$. 

Let $\mathcal{H}$ be the set of Hermitian metrics $H$ such that $r(A_H'', \phi_H) \in \eta_\mu$. Since the neighbourhood $V_\mu$ of Lemma \ref{lem:localstrata} is open and the finite-time gradient flow is continuous in the $C^\infty$ topology by Proposition \ref{lem:continuousdependence}, then $\mathcal{H}$ is open. Let $H^j$ be a sequence of metrics in $\mathcal{H}$ that converge to some Hermitian metric $K$ in the $C^\infty$ topology. The proof of Proposition \ref{lem:convmodgauge} shows that $\| F_A  + [\phi, \phi^*] \|_{L_k^4}$ is bounded along the gradient flow for all $k$, and so Lemma \ref{lem:L4compactness} together with the smooth convergence of $H^j$ shows that there exists a Higgs pair $(A^\infty, \phi^\infty)$, sequences $g_j \in \mathcal{G}$ and $t_j \in \R$ such that $g_j \cdot (A_{H^j}(t_j), \phi_{H^j}(t_j))$ converges to $(A^\infty, \phi^\infty)$ in the $C^\infty$ topology. Let $(A_K^\infty, \phi_K^\infty)$ denote the limit of the gradient flow with initial conditions $(A_K'', \phi_K)$ and note that to prove that $\mathcal{H}$ is closed, it suffices to show that $(A^\infty, \phi^\infty) = (A_K^\infty, \phi_K^\infty)$.

For notation let $H_j = H^j(t_j)$ and $K_j = K(t_j) = H_j h_j$. A calculation shows that
\begin{align}
\begin{split}
(d_{t_j}^{K_j})' - (d_{t_j}^{H_j})' & = h_j^{-1} (d_{t_j}^{H_j})' h_j \\
\phi_{K_j}^* - \phi_{H_j}^* & = h^{-1} [\phi_{H_j}^*, h]
\end{split}
\end{align}
Let $D_j: \Omega^0(\End(E)) \rightarrow \Omega^1(\End(E))$ denote the operator $u \mapsto (d_{t_j}^{H_j})' u + [\phi_{H_j}^*, u]$. The proof of Proposition 6.3 in \cite{Simpson88} shows that the distance measure between metrics $\sup \sigma(H_t, K_t)$ is decreasing with time, and so $\| h_j - \id \|_{C^0} \rightarrow 0$ as $j \rightarrow \infty$. Then we have for any smooth test $1$-form $\beta$
\begin{align}
\begin{split}
\left< (d_{t_j}^{K_j})' - (d_{t_j}^{H_j})' + \phi_{K_j}^* - \phi_{H_j}^* , \beta \right> & = \left< h_j^{-1} \left( (d_{t_j}^{H_j})' h_j + [\phi_{H_j}^*, h_j] \right), \beta \right> \\
& \leq C \left< (d_{t_j}^{H_j})' h_j + [\phi_{H_j}^*, h_j], \beta \right> \\
& =  C \left< h_j, D_j^* \beta \right> \\
& = C \left< h_j, (D_j - D_\infty)^* \beta \right> + C \left< h_j, D_\infty^* \beta \right>
\end{split}
\end{align}
The first term converges to zero since $D_j - D_\infty \rightarrow 0$ smoothly, and since $h_j \rightarrow \id$ in $C^0$ then the second term becomes
\begin{equation}
\left< h_j, D_\infty^* \beta \right> \rightarrow \left< \id, D_\infty^* \beta \right> = \int_X d'^* \tr \beta = 0
\end{equation}
by Stokes' theorem. Therefore $(d_{t_j}^{K_j})' - (d_{t_j}^{H_j})' + \phi_{K_j}^* - \phi_{H_j}^* \rightharpoonup 0$ weakly in $L^2$ and so $(A^\infty, \phi^\infty) = (A_K^\infty, \phi_K^\infty)$. Therefore $\mathcal{H}$ is both open and closed, which completes the proof of Lemma \ref{lem:convergeforallmetrics}.
\end{proof}

\begin{proof}[Proof of Theorem \ref{prop:equivstratification}]

The result of Lemma \ref{lem:convergeforallmetrics} shows that $\mathcal{B}_\mu \subseteq \mathcal{C}_\mu$ for each Harder-Narasimhan type $\mu$. Since the analytic stratification and the algebraic stratification are both partitions of $\mathcal{B}$, then this implies that the two stratifications are equal.
\end{proof}

Next we prove Proposition \ref{lem:stablemanifold}, which provides a description of each stratum in terms of the action of $\mathcal{G}^\C$. Let ${\mathcal{G}_{H^2}^\C}$ denote the completion of the complex gauge group $\mathcal{G}^\C$ in the $H^2$ norm on $\Omega^0(\End(E))$, and note that for $X$ a Riemann surface, the Sobolev embedding theorem shows that ${\mathcal{G}_{H^2}^\C} \subset \mathcal{G}_{C^0}^\C$, the completion of $\mathcal{G}^\C$ in the $C^0$ norm. Let $\mathcal{B}_{H^1}$ denote the completion of the space $\mathcal{B}$ in the $H^1$ norm. For a fixed $C^\infty$ filtration $(*)$, define $\UT(E, *)$ to be the subspace of bundle endomorphisms preserving $(*)$, and similarly let $(\mathcal{G}^\C_{H^2})_*$ denote the subgroup of elements of ${\mathcal{G}_{H^2}^\C}$ which preserve $(*)$. Also let ${(T^* \mathcal{A})_*}$ denote the space of pairs $(A'', \phi)$ such that both $d_A''$ and $\phi$ preserve $(*)$, and let ${(T^* \mathcal{A})_*}_{H^1}$ denote the completion of this space in the $H^1$ norm. Let $\mathcal{B}_*$ and ${\mathcal{B}_*}_{H^1}$ be the respective restrictions of ${(T^* \mathcal{A})_*}$ and ${(T^* \mathcal{A})_*}_{H^1}$ to the space of Higgs pairs. Let $B_*^{ss}$ denote the space of Higgs pairs preserving the filtration $(*)$ such that the pairs $(F_i, \phi_i)$ are semistable for all $i$ with slope strictly decreasing in $i$, where $F_i = E_i / E_{i-1}$ and $\phi_i$ is induced by $\phi$ on $F_i$.

\begin{lem}\label{lem:relategroupfiltration}
If $(*)$ is a filtration of type $\mu$ then $\mathcal{B}_\mu = \mathcal{G}^\C \cdot \mathcal{B}_*^{ss}$ and $(\mathcal{B}_\mu)_{H^1} = \mathcal{G}^\C_{H^2} \cdot (\mathcal{B}_*^{ss})_{H^1}$.
\end{lem}

\begin{proof}
As in the proof of \cite{Daskal92} Lemma 2.10 for holomorphic bundles, we note that $\mathcal{G}^\C \cdot \mathcal{B}_*^{ss} \subseteq \mathcal{B}_\mu$. Conversely, if $(A'', \phi) \in \mathcal{B}_\mu$ then there is a $\phi$-invariant holomorphic filtration of $(E, \phi)$ which is equivalent to $(*)$ by an element of $\mathcal{G}^\C$. The second equality follows in the same way.
\end{proof}

In order to proceed further, we also need the following local description of the space of Higgs bundles close to a point $(A'', \phi) \in \mathcal{B}_{H^1}$. Define the operator 
\begin{equation}\label{eqn:Qinfaction}
\tilde{L}: \Lie(\mathcal{G}^\C_{H^2}) \oplus \Lie(\mathcal{G}^\C_{H^2}) \rightarrow  T_{(A'', \phi)} \left( T^* \mathcal{A} \right)_{H^1}
\end{equation}
by $\tilde{L}(u, v) = \left( \begin{matrix} d_A'' u \\ [\phi, u] \end{matrix} \right) + J \left( \begin{matrix} d_A'' v \\ [\phi, v] \end{matrix} \right) = \rho_\C(u) + J \rho_\C(v)$, where $J$ is the complex structure 
$$J = \left( \begin{matrix} 0 & (\cdot)^* \\ -(\cdot)^* & 0 \end{matrix} \right).$$

Since $\tilde{L}$ is elliptic then $T_{(A'', \phi)} (T^* \mathcal{A})_{H^1} \cong \im \tilde{L} \oplus \ker \tilde{L}^*$. The following lemma shows that when $(A'', \phi) \in \mathcal{B}_{H^1}$ then the same is true for the operator $\rho_\C : \Lie(\mathcal{G}^\C_{H^2}) \rightarrow T_{(A'', \phi)} \left( T^* \mathcal{A} \right)_{H^1}$.

\begin{lemma}
Let $(A_0'', \phi_0) \in \mathcal{B}_{H^1}$. Then $T_{(A_0'', \phi_0)} (T^* \mathcal{A})_{H^1} = \im \rho_\C \oplus \ker \rho_\C^*$.
\end{lemma}

\begin{proof}
Since $\tilde{L}$ is elliptic then $\im \tilde{L} = \im \rho_\C + \im J \rho_\C$ is closed, and so we have $T_{(A_0'', \phi_0)} (T^* \mathcal{A})_{H^1} = \im \tilde{L} \oplus \ker \tilde{L}^*$. The adjoint $\tilde{L}^*$ is given by $\tilde{L}^*(a'', \varphi) = (\rho_\C^*(a'', \varphi), -\rho_\C^* J (a'', \varphi))$ and so $\ker \tilde{L}^* = \ker \rho_\C^* \cap \ker (\rho_\C^* J)$. 

Since $\im \rho_\C + \im J \rho_\C$ is closed and $(A'', \phi) \in \mathcal{B}_{H^1}$, then $\im \rho_\C \subseteq ( \im J \rho_\C )^\perp$ and $\im J \rho_\C \subseteq (\im \rho_\C)^\perp$. Lemma \ref{lem:closedsubspace} then shows that  $\im \rho_\C = ( \im J \rho_\C )^\perp$ and $\im J \rho_\C = (\im \rho_\C)^\perp$, so $\im \rho_\C$ and $\im J \rho_\C$ are closed subspaces of $\im \tilde{L}$ and we have $\im \tilde{L} = \im \rho_\C \oplus \im J \rho_\C$. Therefore
$$
T_{(A_0'', \phi_0)} (T^* \mathcal{A})_{H^1} = \im \rho_\C \oplus \im J \rho_\C \oplus \left( \ker \rho_\C^* \cap \ker (\rho_\C^* J) \right)
$$
Since $\im J \rho_\C \oplus \left( \ker \rho_\C^* \cap \ker (\rho_\C^* J) \right) \subseteq \ker \rho_\C^* \subseteq (\im \rho_\C)^\perp$ then applying Lemma \ref{lem:closedsubspace} again shows that the set inclusions are in fact equalities, which gives the decomposition $T_{(A_0'', \phi_0)} (T^* \mathcal{A})_{H^1} = \im \rho_\C \oplus \ker \rho_\C^*$.
\end{proof}

The next lemma follows from the inverse function theorem.

\begin{lemma}\label{lem:inversefunction}
The map $f : (\ker \rho_\C)^\perp \times \ker \rho_\C^* \rightarrow T^* \mathcal{A}_{H^1}$ given by $f(u, a'', \varphi) = e^u \cdot (A'' + a'', \phi + \varphi)$ is a local diffeomorphism at $(0,0,0)$.
\end{lemma}

\begin{proof}
The derivative of $f$ at $(0,0,0)$ is the map $df (\delta u, \delta a'', \delta \varphi) = \rho_\C(\delta u) + (\delta a'', \delta \varphi)$, which is an isomorphism by the previous lemma. The inverse function theorem then shows that $f$ is a local diffeomorphism.
\end{proof}

Now let $\mathcal{S}_{(A'', \phi)}$ be the slice given by
$$
\mathcal{S}_{(A'', \phi)} = \ker \rho_\C^* \cap \left\{ (a'', \varphi) \in T_{(A'', \phi)} (T^* \mathcal{A})_{H^1} \, : \, d_A'' \varphi + [a'', \phi] + [a'', \varphi] = 0 \right\}
$$

\begin{lemma}
Let $\tilde{f}$ be the restriction of $f$ to $(\ker \rho_\C)^\perp \times \mathcal{S}_{(A'', \phi)}$. If $(\tilde{A}'', \tilde{\phi}) \in \mathcal{B}_{H^1}$ satisfies $\| (\tilde{A}'', \tilde{\phi}) - (A'', \phi) \|_{H^1} < \varepsilon$ then there exist unique elements $u \in (\ker \rho_\C)^\perp$ and $(a'', \varphi) \in \mathcal{S}_{(A'', \phi)}$ such that $(A'', \phi) = \tilde{f}(u, a'', \varphi)$.
\end{lemma}

\begin{proof}
Lemma \ref{lem:inversefunction} shows that there exists $(u, a'', \varphi) \in (\ker \rho_\C)^\perp \times \ker \rho_\C^*$ such that $(\tilde{A}'', \tilde{\phi}) = f(u, a'', \varphi)$. Therefore only remains to show that $(a'', \varphi) \in \mathcal{S}_{(A'', \phi)}$, which results from observing that $(\tilde{A}'', \tilde{\phi}) \in \mathcal{B}_{H^1}$ iff $e^u \cdot (A'' + a'', \phi + \varphi) \in \mathcal{B}_{H^1}$ iff $(A'' + a'', \phi + \varphi) \in \mathcal{B}_{H^1}$.
\end{proof}

\begin{prop}\label{prop:localhomeo}
Fix $(A'', \phi) \in \mathcal{B}_{H^1}$. Then the map $\tilde{f} : (\ker \rho_\C)^\perp \times \mathcal{S}_{(A'', \phi)} \rightarrow \mathcal{B}_{H^1}$ is a local homeomorphism from a neighbourhood of zero in $(\ker \rho_\C)^\perp \times \mathcal{S}_{(A'', \phi)}$ to a neighbourhood of $(A'', \phi) \in \mathcal{B}_{H^1}$.
\end{prop}

\begin{proof}[Proof of Proposition \ref{prop:localhomeo}]
If $(a'', \varphi) \in \mathcal{S}_{(A'', \phi)}$ then $f(u, a'', \varphi) \in \mathcal{B}_{H^1}$ for any $u \in (\ker \rho_\C^*)^\perp$, which combined with the previous lemma shows that $\tilde{f}$ is surjective onto a neighbourhood of $(A'', \phi) \in \mathcal{B}_{H^1}$. Since $\tilde{f}$ is the restriction of a local diffeomorphism then it is a local homeomorphism onto a neighbourhood of $(A'', \phi)$ in $\mathcal{B}_{H^1}$.
\end{proof}

Given a filtration $(*)$ of the bundle $E$, define the subset of the slice consisting of variations that preserve the filtration by $\left( \mathcal{S}_{(A'', \phi)}\right)_* =  \mathcal{S}_{(A'', \phi)} \cap \Omega^{0,1}(\UT(E, *)) \oplus \Omega^{1,0}(\UT(E, *))$. Let $p$ be the projection $p : \left( \ker \rho_\C \right)^\perp \times \mathcal{S}_{(A'', \phi)}   \rightarrow \left( \ker \rho_\C \right)^\perp \times  \left( \mathcal{S}_{(A'', \phi)} \right)_*$. We then have the following description of each stratum close to a critical point.

\begin{lem}\label{lem:cutoutstrata}
Let $(A_0, \phi_0) \in \mathcal{B}_{H^1}$ be a critical point of $\YMH$ with Harder-Narasimhan filtration $(*)$. Then there exists $\varepsilon > 0$ such that for any $(A'', \phi) \in (\mathcal{B}_\mu)_{H^1}$ with $\hone{ (A'', \phi) - (A_0'', \phi_0) } < \varepsilon$, there exists $\left( u, a'', \varphi \right) \in \ker(1-p)$ such that $\tilde{f}(u,a'', \varphi) = (A'', \phi)$.
\end{lem}

\begin{rek}
Conversely, this lemma implies that if $(A'', \phi) \in \mathcal{B}_{H^1} \setminus (\mathcal{B}_\mu)_{H^1}$ and $(A'', \phi)$ is close to $(\mathcal{B}_\mu)_{H^1}$, then there exists $\left( u, a'', \varphi  \right)$ satisfying $(1-p)(a'', \varphi) \neq 0$ and $\tilde{f}\left( u, a'', \varphi \right) = (A'', \phi)$. In other words we have a criterion that describes exactly when a point in a neighbourhood of $(A_0'', \phi_0)$ lies in the stratum $(\mathcal{B}_\mu)_{H^1}$.
\end{rek}

\begin{proof}[Proof of Lemma \ref{lem:cutoutstrata}]

Proposition \ref{prop:localhomeo} states that there exists $\varepsilon > 0$ such that given a point $(A'', \phi)$ within a distance $\varepsilon$ from $(A_0'', \phi_0)$ in the $H^1$ norm there exists  $(u, a'', \varphi) \in \left( \ker \rho_\C \right)^\perp \oplus \mathcal{S}_{(A'', \phi)}$ such that $\tilde{f}(u, a'', \varphi) = (A'', \phi)$. Restricting to the stratum $\mathcal{B}_\mu$ we follow the same steps as in the proof of Proposition 3.5 from \cite{Daskal92} (for the Yang-Mills functional and unitary connections), except for the functional $\YMH$ and $GL(n, \C)$ connections, to show that $(a'', \varphi) \in \left( \mathcal{S}_{(A'', \phi)} \right)_*$. Therefore the projection $p$ is the identity on this space, which completes the proof.
\end{proof}

Let $\left( \ker \rho_\C \right)_*^\perp = \left( \ker \rho_\C \right)^\perp \cap \Omega^0 (\UT(E, *))$. The previous lemma describes a neighbourhood in $(\mathcal{B}_\mu)_{H^1}$, and now we describe a neighbourhood in $(\mathcal{B}_*^{ss})_{H^1}$.

\begin{lemma}\label{lem:homeofixedfiltration}
The restricted map $\tilde{f}_* : \left( \ker \rho_\C \right)_*^\perp \times  \left( \mathcal{S}_{(A'', \phi)} \right)_* \rightarrow (\mathcal{B}_*^{ss})_{H^1}$ is a local homeomorphism.
\end{lemma}

\begin{proof}
Clearly $\tilde{f}_*$ maps into $(\mathcal{B}_*^{ss})_{H^1}$. Since it is the restriction of a local homeomorphism then it is a local homeomorphism onto its image, and so the proof reduces to showing that $\tilde{f}_*$ is locally surjective. Lemma \ref{lem:cutoutstrata} shows that if $(\tilde{A}'', \tilde{\phi})$ is close to $(A'', \phi)$ in the $H^1$ norm then there exists $u \in \left( \ker \rho_\C \right)^\perp$ and $(a'', \varphi) \in \left( \mathcal{S}_{(A'', \phi)} \right)_*$ such that $e^u \cdot (A'' + a'', \phi + \varphi) = (\tilde{A}'', \tilde{\phi})$. The proof then reduces to showing that $u \in \left( \ker \rho_\C \right)_*^\perp$. Restricting our viewpoint to the holomorphic structures, we see that a weak sub-bundle $\pi$ corresponding to a term in the Harder-Narasimhan filtration $(*)$ is holomorphic, and so the equation in Lemma (3.2) of \cite{Daskal92} holds for $\pi$. This allows us to prove a Higgs-bundle version of Lemma (3.3) in \cite{Daskal92}, which shows that $u \in \left( \ker \rho_\C \right)_*^\perp$.
\end{proof}

\begin{prop}\label{lem:stablemanifold}
$\left( \mathcal{B}_\mu \right)_{H^1}$ is homeomorphic to
\begin{equation*}
\mathcal{G}_{H^2}^{\C} \times_{(\mathcal{G}^{\C}_{H^2})_*} \left( \mathcal{B}_*^{ss} \right)_{H^1} \cong \mathcal{G}_{H^2} \times_{{\mathcal{G}_{diag}}_{H^2}} \left( \mathcal{B}_*^{ss} \right)_{H^1}
\end{equation*}
where $\mathcal{G}_{diag} \subset \mathcal{G}$ denotes the space of diagonal gauge transformations with respect to the fixed $C^\infty$ filtration $(*)$.
\end{prop}

\begin{proof}[Proof of Proposition \ref{lem:stablemanifold}]

Define the map $\psi: \mathcal{G}^\C_{H^2} \times_{(\mathcal{G}^\C_{H^2})_*} {\mathcal{B}_*^{ss}}_{H^1} \rightarrow (\mathcal{B}_\mu)_{H^1}$ by $\psi([g, (A'', \phi)]) = g \cdot (A'', \phi)$. If $\psi([g_1, (A_1'', \phi_1)]) = \psi([g_2, (A_2'', \phi_2)])$ then $g_1 \cdot (A_1'', \phi_1) = g_2 \cdot (A_2'', \phi_2)$ with $(A_1'', \phi_1), (A_2'', \phi_2) \in (\mathcal{B}_*^{ss})_{H^1}$, so $g_1^{-1} g_2 \in {\mathcal{G}_*^{\C}}_{H^2}$ and therefore $\psi$ is injective. Lemma \ref{lem:relategroupfiltration} shows that $\psi$ is surjective onto $(\mathcal{B}_\mu)_{H^1}$, and so the first equality in the proposition will follow if we can show that $\psi$ is a local homeomorphism.

Lemma \ref{lem:homeofixedfiltration} shows that a neighbourhood of a point $(A_0'', \phi_0) \in (\mathcal{B}_*^{ss})_{H^1}$ is homeomorphic to a neighbourhood of zero in $\left( \ker \rho_\C \right)_*^\perp \times \left( \mathcal{S}_{(A_0'', \phi_0)} \right)_*$. Therefore, in $\mathcal{G}_{H^2}^{\C} \times_{{\mathcal{G}_*^{\C}}_{H^2}} \left( \mathcal{B}_*^{ss} \right)_{H^1}$ we have $\left[ g, (A'', \phi) \right] = \left[ g, e^u \cdot (A_0'' + a'', \phi_0 + \varphi) \right] = \left[ e^{-u} g, A_0'' + a'', \phi_0 + \varphi \right]$, since $e^u \in (\mathcal{G}^{\C}_{H^2})_*$. This implies that $\psi(\left[ g, (A'', \phi) \right]) = e^{-u} g \cdot (A_0'' + a'', \phi_0 + \varphi)$ with $(a'', \varphi) \in \left( \mathcal{S}_{(A_0'', \phi_0)} \right)_*$. Lemma \ref{lem:cutoutstrata} then shows that $\psi$ is a local homeomorphism when $g$ is close to the identity, and translating this result by the action of the complex gauge group shows that $\psi$ is a local homeomorphism for all $g$.

The homeomorphism $\mathcal{G}_{H^2}^\C \cong \left( \mathcal{G}_*^{\C} \right)_{H^2} \times_{{\mathcal{G}_{diag}}_{H^2}} \mathcal{G}_{H^2}$ from Theorem 2.16 in \cite{Daskal92} completes the proof of the second equality in the statement of Proposition \ref{lem:stablemanifold}.
\end{proof}

\begin{corollary}
$$\mathcal{B}_\mu \cong \mathcal{G}^\C \times_{\mathcal{G}_*^\C} \mathcal{B}_*^{ss} \cong \mathcal{G} \times_{\mathcal{G}_{diag}} \mathcal{B}_*^{ss} $$
\end{corollary}

\begin{proof}
Lemma 14.8 of \cite{AtiyahBott83} shows that every $\mathcal{G}_{H^2}^\C$-orbit in $\mathcal{A}_{H^1}^{0,1}$ contains a $C^\infty$ holomorphic structure $d_A''$. If the holomorphic structure $A''$ is smooth, then since the Higgs bundle equation $d_A'' \phi = 0$ is elliptic then all $\phi$ satisfying this condition are smooth. Therefore every $\mathcal{G}_{H^2}^\C$-orbit in $\mathcal{B}_{H^1}$ contains a $C^\infty$ Higgs pair. Moreover, if two $C^\infty$ holomorphic structures $A_1''$ and $A_2''$ are isomorphic by an element $g \in \mathcal{G}_{H^2}$, then bootstrapping the equation $g A_1'' - A_2'' g = d'' g$ shows that $g$ is smooth also, and so every $\mathcal{G}_{H^2}^\C$-orbit in $\mathcal{B}_{H^1}$ contains exactly one $\mathcal{G}^\C$ orbit of smooth Higgs pairs. The corollary then follows from Proposition \ref{lem:stablemanifold}.
\end{proof}

\section{Convergence to the Graded Object of the filtration}\label{sec:gradedobject}

The results of Section \ref{sec:analysis} show that the gradient flow of $\YMH$ converges smoothly to a critical point of $\YMH$, and the results of Section \ref{sec:stratification} describe the type of the $\phi$-invariant Harder-Narasimhan filtration at the limit. The purpose of this section is to provide an algebraic description of the isomorphism class of the limit of the gradient flow, a Higgs bundle version of Corollary (5.19) of \cite{Daskal92} (for the Yang-Mills functional on a Riemann surface) and Theorem 1 of \cite{DaskalWentworth04} (Yang-Mills on a K\"ahler surface). To describe the limit algebraically requires a description of the appropriate Higgs bundle versions of the Seshadri filtration and the Harder-Narasimhan-Seshadri filtration, which is contained in the following Propositions (cf \cite{DaskalWentworth04} Propositions 2.5 and 2.6 for holomorphic bundles)
\begin{prop}\label{prop:seshadrifiltration}
Let $(A'', \phi)$ be a Higgs-semistable structure on $E$. Then there is a filtration of $E$ by $\phi$-invariant holomorphic sub-bundles
\begin{equation*}
0 = F_0 \subset F_1 \subset \cdots \subset F_\ell = E
\end{equation*}
called a \emph{$\phi$-invariant Seshadri filtration} of $E$, such that $F_i / F_{i-1}$ is Higgs stable for all $i$ (with respect to the Higgs structure induced from $(A'', \phi)$), and $\mu(F_i / F_{i-1}) = \mu(E)$. The graded object $\displaystyle{\Gr^{\mathrm{S}}(A'', \phi) = \bigoplus_{i=1}^\ell F_i / F_{i-1}}$ is uniquely determined by the isomorphism class of $(A'', \phi)$.
\end{prop}

\begin{prop}\label{prop:hnsfiltration}
Let $(A'', \phi)$ be a Higgs structure on $E$. Then there is a double filtration $\{ E_{i,j} \}$ of $E$, called a \emph{$\phi$-invariant Harder-Narasimhan-Seshadri filtration} of $E$ (abbr. HNS filtration) such that if $\{ E_i \}_{i=1}^\ell$ is the $\phi$-invariant HN filtration of $E$ then
\begin{equation*}
E_{i-1} = E_{i,0} \subset E_{i,1} \subset \cdots \subset E_{i, \ell_i} = E_i
\end{equation*}
is a Seshadri filtration of $E_i / E_{i-1}$. The associated graded object
\begin{equation}\label{eqn:defofgradedobject}
\Gr^{\HNS} (A'', \phi) = \bigoplus_{i=1}^\ell \bigoplus_{j=1}^{\ell_i} Q_{i,j}
\end{equation}
is uniquely determined by the isomorphism class of $(A'', \phi)$.
\end{prop}

Recall the gradient flow retraction $r: \mathcal{B} \rightarrow \mathcal{B}_{crit}$ onto the set of critical points $\mathcal{B}_{crit}$ defined in Theorem \ref{thm:convergence}. The main theorem of this section is the following
\begin{thm}\label{thm:gradedobject}
The isomorphism class of the gradient flow retraction is given by
\begin{equation}\label{eqn:retractiongradedobject}
r(A'', \phi) \cong \Gr^{\HNS}(A'', \phi)
\end{equation}
\end{thm}

Consider a sequence $t_n \rightarrow \infty$, and denote $(A(t_n)'', \phi(t_n))$ by $(A_n'', \phi_n)$. Let $g_n \in \mathcal{G}^\C$ be the complex gauge transformation corresponding to the finite-time gradient flow from  time $t_0$ to $t_n$, i.e. $(A_n'', \phi_n) = g_n \cdot (A_0'', \phi_0)$. Let $S$ be the first term in the Harder-Narasimhan-Seshadri filtration of $E$, and let $f_0: S \hookrightarrow E$ be the $\phi$-invariant, holomorphic inclusion. Define the map $f_n : S \hookrightarrow E$ by $f_n = g_n \circ f_0$, and note that since $f_0$ and $g_n$ are $\phi$-invariant holomorphic sections of the associated Higgs bundles $\Hom(S, E_0)$ and $\Hom(E_0, E_n)$ (with the induced Higgs fields) then $f_n$ is also holomorphic and $\phi$-invariant. Define the operators 
\begin{align*}
& D_n'' : \Omega^0 (\ad(E)) \rightarrow \Omega^{0, 1}(\End(E)) \oplus \Omega^{1,0}(\End(E)) \\
\mathrm{and} \quad  & D_{i,j} : \Omega^0(\Hom(E_i, E_j) \rightarrow \Omega^{0,1} \left(\Hom(E_i, E_j) \right) \oplus \Omega^{1,0} \left(\Hom(E_i, E_j) \right)
\end{align*}
by $u \mapsto \left(  d_{A_n}'' u, [\phi_n, u] \right)$, and $u \mapsto \left( d''u + A_j'' u - u A_i'', \phi_j u - u \phi_i \right)$. Let $g_{i,j} \in \mathcal{G}^\C$ correspond to the finite-time gradient flow from time $t_i$ to $t_j$ (i.e. $g_{i,j} \cdot (A_i, \phi_i) = (A_j, \phi_j)$). Then a simple calculation shows that $D_{i,j} g_{i,j} = 0$. The proof of Proposition \ref{lem:convmodgauge} shows that $\| \phi \|_{H^k}$ and $\| F_A \|_{H^k}$ are bounded for all $k$ along the gradient flow of $\YMH$, and so for all $\ell$
\begin{equation}\label{eqn:ellipticestimate}
\| D_{i,j} u \|_{H^{\ell-1}} \leq C \| u \|_{H^\ell},
\end{equation}
where the bound $C$ is uniform along the gradient flow, by \cite{Rade92} Proposition A and Lemma \ref{lem:L4compactness} in this paper. After these preliminaries we can now prove the following claim.

\begin{cla}
$f_n$ converges in the $H^k$ norm for all $k$ to some non-zero $\phi$-invariant holomorphic map $f_\infty$.
\end{cla}

\begin{proof}

Replace $f_n$ by $\frac{f_n}{\ltwo{f_n}}$ (note that $\ltwo{f_n} \neq 0$ for all $n$ since $\| f_0 \|_{L^2} \neq 0$ and $g_n$ is an automorphism of $E$) and consider $D_{0,n}'' f_n = D_{0, \infty}'' f_n + [\beta_n, f_n]$ where $\beta_n \rightarrow 0$ in $H^k$ for all $k$ (since $D_n'' \rightarrow D_\infty''$ in $H^k$ for all $k$). Since $f_n$ is holomorphic, then $D_{0,n}'' f_n = 0$. Therefore for any $\ell$ we have the estimate
\begin{equation}\label{eqn:derivativebound}
\left\| D_{0, \infty}'' f_n \right\|_{H^\ell} \leq \left\| \beta_n \right\|_{C^0} \left\| f_n \right\|_{H^\ell}
\end{equation}
Since $\beta_n \rightarrow 0$ smoothly then along a subsequence (also denoted $f_n$), $f_n$ bounded in $H^\ell$ implies that $f_n$ is bounded in $H^{\ell + 1}$, where the bound only depends on $\| f_n \|_{H^\ell}$. Since $\| f_n \|_{L^2} = 1$, then by induction $\| f_n \|_{H^\ell} \leq C_\ell$ for all $\ell$. Therefore there exists $f_\infty$ such that $f_n \rightarrow f_\infty$ strongly in $H^{\ell-1}$ for all $\ell$. The estimate \eqref{eqn:ellipticestimate} for the operator $D_{0, \infty}''$ shows that
\begin{align*}
\left\| D_{0, \infty}'' f_\infty \right\|_{H^{\ell-1}} & \leq \left\| D_{0, \infty}'' (f_n - f_\infty) \right\|_{H^{\ell-1}} + \left\| D_{0, \infty}'' f_n \right\|_{H^{\ell-1}} \\
 & \leq C \left\| f_n - f_\infty \right\|_{H^\ell} + \left\| D_{0, \infty}'' f_n \right\|_{H^{\ell-1}}
\end{align*}
Since $\beta_n \rightarrow 0$ and $\| f_n \|_{H^\ell}$ is bounded, then \eqref{eqn:derivativebound} shows that the right-hand side of the above estimate approaches zero as $n \rightarrow \infty$ for all $\ell$. Therefore $D_{0, \infty}'' f_\infty = 0$ and so $f_\infty$ is holomorphic. Since $\| f_n \|_{L^2} = 1$ for all $n$ then $f_\infty \neq 0$.
\end{proof}

Theorem \ref{prop:equivstratification} shows that the \emph{type} of the Harder-Narasimhan filtration is preserved in the limit. The next result shows that the destabilising Higgs sub-bundles in the Harder-Narasimhan filtration along the gradient flow converge to the destabilising Higgs sub-bundles of the limiting Higgs pair. In the following we use the projection $\pi : E \rightarrow E$ to denote the sub-bundle $\pi(E)$.

\begin{prop}\label{prop:HNbundlepreserved}
Let $\{ \pi_t^{(i)} \}$ be the HN filtration of a solution $(A_t'', \phi_t)$ to the gradient flow equations \eqref{eqn:gradfloweqns}, and let $\{ \pi_\infty^{(i)} \}$ be the HN filtration of the limit $(A_\infty'', \phi_\infty)$. Then there exists a subsequence $\{ t_j \}$ such that $\pi_{t_j}^{(i)} \rightarrow \pi_\infty^{(i)}$ in $L^2$ for all $i$.
\end{prop}

To prove this we need the following lemmas.

\begin{lemma}
$\| D_t''(\pi_t^{(i)}) \|_{L^2} \rightarrow 0$
\end{lemma}

\begin{proof}
Let $D_{t_j}'' : \Omega^0(\End(E)) \rightarrow \Omega^{0,1}(\End(E)) \oplus \Omega^{1,0}(\End(E))$ denote the infinitesimal action of $\mathcal{G}^\C$ at time $t$, i.e. $D_{t}''(u) = ( d_{A_t}'' u, [\phi_t, u])$. The Chern-Weil formula of \cite{Simpson88} shows that 
\begin{equation}\label{eqn:simpsonchernweil}
\deg(\pi_{t}^{(i)}) = \frac{\sqrt{-1}}{2\pi} \int_X \tr \left( \pi_t^{(i)} * (F_{A_t} + [\phi_t, \phi_t^*]) \right) - \| D_t''(\pi_t^{(i)}) \|_{L^2}^2
\end{equation}
Along the finite-time flow $d_i = \deg(\pi_{t}^{(i)})$ is fixed, therefore we can re-write \eqref{eqn:simpsonchernweil}
\begin{multline}\label{eqn:simpsonchernweilexpanded}
\| D_t''(\pi_t^{(i)}) \|_{L^2} = -d_i + \frac{\sqrt{-1}}{2 \pi} \int_X \tr \left( \pi_t^{(i)} *(F_{A_\infty} + [\phi_\infty, \phi_\infty^*] ) \right) \\ + \frac{\sqrt{-1}}{2\pi} \int_X \tr \left( \pi_t^{(i)} * \left(F_{A_t} + [\phi_t, \phi_t^*] - F_{A_\infty} - [\phi_\infty, \phi_\infty^*] \right) \right)
\end{multline}
Theorem \ref{thm:convergence} shows that $F_{A_t} + [\phi_t, \phi_t^*] \rightarrow F_{A_\infty} + [\phi_\infty, \phi_\infty^*]$ in the $C^\infty$ topology, and therefore since $\pi_t^{(i)}$ is uniformly bounded in $L^2$ (it is a projection) then the last term in \eqref{eqn:simpsonchernweilexpanded} converges to zero. Let $\mu$ be the HN type of $(A_\infty'', \phi_\infty)$. Since $(A_\infty'', \phi_\infty)$ is a critical point of $\YMH$ then we also have 
\begin{equation}
\frac{\sqrt{-1}}{2 \pi} \int_X \tr \left( \pi_t^{(i)} *(F_{A_\infty} + [\phi_\infty, \phi_\infty^*] ) \right) \leq \sum_{k \leq \rank(\pi_\infty^{(i)})} \mu_k = d_i
\end{equation}
Combining all of these results, we see that $\| D_t''(\pi_t^{(i)}) \|_{L^2} \rightarrow 0$.
\end{proof}

In particular, this lemma shows that $\| \pi_{t_j}^{(i)} \|_{H^1} \leq C$ and so there exists some $\tilde{\pi}_\infty^{(i)}$ and a subsequence $t_j$ such that $\pi_{t_j}^{(i)} \rightarrow \tilde{\pi}_\infty^{(i)}$ weakly in $H^1$ and strongly in $L^2$.

\begin{lemma}
$\| D_\infty'' (\tilde{\pi}_\infty^{(i)}) \|_{L^2} = 0$
\end{lemma}

\begin{proof}
$ \| D_\infty'' (\pi_{t_j}^{(i)}) \|_{L^2} \leq \| D_\infty'' (\pi_{t_j}^{(i)}) - D_{t_j}'' (\tilde{\pi}_{t_j}^{(i)}) \|_{L^2} + \| D_{t_j}'' (\tilde{\pi}_{t_j}^{(i)}) \|_{L^2}$. Theorem \ref{thm:convergence} and the previous lemma then show that $\| D_\infty'' (\pi_{t_j}^{(i)}) \|_{L^2} \rightarrow 0$. Since $\pi_{t_j}^{(i)} \rightarrow \tilde{\pi}_\infty^{(i)}$ weakly in $H^1$ then $\| D_\infty'' (\tilde{\pi}_\infty^{(i)}) \|_{L^2} = 0$.
\end{proof}

\begin{lemma}
$\deg(\tilde{\pi}_\infty^{(i)}) = \deg(\pi_\infty^{(i)})$
\end{lemma}

\begin{proof}
The previous lemma and equation \eqref{eqn:simpsonchernweil} show that 
\begin{align}
\begin{split}
\deg(\tilde{\pi}_\infty^{(i)}) & =  \frac{\sqrt{-1}}{2\pi} \int_X \tr \left( \tilde{\pi}_\infty^{(i)} * (F_{A_\infty} + [\phi_\infty, \phi_\infty^*]) \right) \\
  & = \lim_{j \rightarrow \infty} \| D_{t_j}'' \pi_{t_j}^{(i)} \|_{L^2}^2 + \deg( \pi_{t_j}^{(i)} ) \\
  & = \deg( \pi_\infty^{(i)} )
\end{split}
\end{align}
where in the last step we use the result of Theorem \ref{prop:equivstratification} that the type of the HN filtration is preserved in the limit.
\end{proof}

\begin{proof}[Proof of Proposition \ref{prop:HNbundlepreserved}]

The results of the preceding lemmas show that the degree and rank of $\pi_\infty^{(i)}$ and $\tilde{\pi}_\infty^{(i)}$ are the same. For $i=1$, $\pi_\infty^{(1)}$ is the maximal destabilising Higgs sub-bundle of $(A_\infty, \phi_\infty)$, which is the unique Higgs sub-bundle of this degree and rank. Therefore $\pi_\infty^{(1)} = \tilde{\pi}_\infty^{(1)}$. Proceeding by induction on the HN filtration as in \cite{DaskalWentworth04} completes the proof of Proposition \ref{prop:HNbundlepreserved}.
\end{proof}

Following the idea in part (2) of the proof of Lemma 4.5 in \cite{DaskalWentworth04} in the Yang-Mills case, we see that the same argument applies to the Seshadri filtration of a semistable Higgs bundle, except that because of the lack of uniqueness of the Seshadri filtration we can only conclude that the degree and rank of the limiting sub-bundle are the same.

The following lemma is completely analogous to the proof of  (V.7.11) in \cite{Kobayashi87} for holomorphic bundles and so the proof is omitted.

\begin{lem}
Let $(S_1, \phi_1)$  be a stable Higgs bundle, and let $(S_2, \phi_2)$ be a semistable Higgs bundle over a compact Riemann surface $X$. Also suppose that $\frac{\deg(S_1)}{\rank(S_1)} = \frac{\deg(S_2)}{\rank(S_2)}$, and let $f: S_1 \rightarrow S_2$ be a holomorphic map satsifying $f \circ \phi_1 = \phi_2 \circ f$. Then either $f=0$ or $f$ is injective.
\end{lem}

Since the Harder-Narasimhan filtration is preserved in the limit then $(S, A_0, \phi_0)$ is Higgs-stable and $(S, A_\infty, \phi_\infty)$ is Higgs-semistable with the same degree/rank ratio, so the non-zero map $f_\infty$ must be injective. Therefore $\im f_\infty = (S, A_\infty, \phi_\infty)$ is Higgs-stable. Using \cite{Daskal92} Lemma 5.12 we can assume (after unitary co-ordinate changes) that the operator $D_i''$ preserves the bundle $S_\infty$ for all $i$.  To complete the induction we need the following result for the quotient bundle $Q$.

\begin{cla}
Let $Q_k = E_k/S_k$. Then $Q_k = h_k \cdot Q_0$ for some $h_k \in \mathcal{G}^\C(Q)$, the induced connections ${D_j''}^Q$ converge to some ${D_\infty''}^Q$ in the $C^\infty$ norm, and $Q_0$ and $Q_\infty$ have the same $\phi$-invariant Harder-Narasimhan type.
\end{cla}

\begin{proof}
The construction of $h_k$ follows from the following commutative diagram
\begin{equation}
\xymatrix{
0 \ar[r] & S_0 \ar[r] \ar[d]^{f_k} & E_0 \ar[r] \ar[d]^{g_k} & Q_0 \ar@{-->}[d]^{h_k} \ar[r] & 0 \\
0 \ar[r] & S_k \ar[r] & E_k \ar[r] & Q_k \ar[r] & 0
}
\end{equation}
where the map $h_k$ is constructed from the maps $f_k$ and $g_k$ using the exactness of the rows in the diagram.

Using the notation from Lemma 5.12 of \cite{Daskal92}, the induced connection on $Q_k$ is given by ${D_k''}^Q = \tilde{\pi_k}^\perp D_k'' \tilde{\pi_k}^\perp$. Lemma 5.12 from \cite{Daskal92} states that $\tilde{\pi_k} = \pi_\infty$ is constant with respect to $k$, and so ${D_k''}^Q = \pi_\infty^\perp D_k'' \pi_\infty^\perp$ converges to $\pi_\infty^\perp D_\infty'' \pi_\infty^\perp = {D_\infty''}^Q$. Finally, Theorem \ref{prop:equivstratification} shows that $Q_0$ and $Q_\infty$ have the same Harder-Narasimhan type.
\end{proof}

Therefore we can apply the previous argument to the first term in the double filtration of $Q$. Repeating this process inductively shows that the limit of the gradient flow $\YMH$ along the sequence $\{ t_n \}$ is the graded object associated to the $\phi$-invariant Harder-Narasimhan-Seshadri filtration of $(A'', \phi)$. Since Theorem \ref{thm:convergence} shows that the limit exists along the flow independently of the subsequence chosen then the limit is $\Gr^{\HNS}(A'', \phi)$, completing the proof of Theorem \ref{thm:gradedobject}.



\def\cprime{$'$}

\end{document}